\documentclass[aos,MSNbibl,seceqn,citesort,dvips]{arximspdf}
\usepackage{dcolumn}
\usepackage{graphicx}

\doi{10.1214/11-AOS947}
\volume{40}
\issue{1}
\pubyear{2012}
\firstpage{73}
\lastpage{103}

\makeatletter

\newcolumntype{d}[1]{D{.}{.}{#1}}

\newtheorem{theorem}{Theorem}[section]
\newtheorem{lemma}{Lemma}[section]

\newproclaim{definition}{Definition}[section]
\newproclaim{Experiment}{Experiment}

\newcommand{\eps}{\varepsilon}
\newcommand{\goto}{\rightarrow}
\newcommand{\bphi}{\bar{\Phi}}
\newcommand{\sgn}{\operatorname{sgn}}
\newcommand{\hb}{\hat{\beta}}
\newcommand{\hamm}{\operatorname{Hamm}}
\newcommand{\call}{{\mathcal I}_0}
\newcommand{\cf}{\bar{F}}
\newcommand{\ty}{\tilde{Y}}

\setattribute{abstract}{width}{300pt}

\makeatother

\begin{document}
\begin{frontmatter}

\title{UPS delivers optimal phase diagram in high-dimensional variable selection}
\runtitle{UPS delivers optimal phase diagram}

\begin{aug}
\author[A]{\fnms{Pengsheng} \snm{Ji}\thanksref{t1,t2}\ead[label=e1]{pj54@cornell.edu}}
\and
\author[B]{\fnms{Jiashun} \snm{Jin}\corref{}\thanksref{t1}\ead[label=e2]{jiashun@stat.cmu.edu}}
\runauthor{P. Ji and J. Jin}
\affiliation{Cornell University and Carnegie Mellon University}
\address[A]{Department of Statistical Science\\
Cornell University\\
Ithaca, New York 14853\\
USA\\
\printead{e1}} 
\address[B]{Department of Statistics\\
Carnegie Mellon University\\
Pittsburgh, Pennsylvania 15213 \\
USA\\
\printead{e2}}
\end{aug}

\thankstext{t1}{Supported in part by NSF CAREER Award DMS-09-08613.}

\thankstext{t2}{Supported in part by NSF Grant DMS-08-05632.}

\received{\smonth{5} \syear{2011}}
\revised{\smonth{11} \syear{2011}}

%
\begin{abstract}
Consider a linear model $Y = X \beta+ z$, $z \sim N(0, I_n)$. Here, $X
= X_{n, p}$, where both $p$ and $n$ are large, but $p > n$. We model
the rows of~$X$ as i.i.d. samples from $N(0, \frac{1}{n} \Omega)$,
where $\Omega$ is a $p\times p$ correlation matrix, which is unknown to
us but is presumably sparse. The vector~$\beta$ is also unknown but has
relatively few nonzero coordinates, and we are interested in
identifying these nonzeros.

We propose the Univariate Penalization Screeing (UPS)
for variable selection. This is a screen and clean method where we
screen with univariate thresholding and clean with penalized MLE. It
has two important properties: sure screening and separable after
screening. These properties enable us to reduce the original
regression problem to many small-size regression problems that can be
fitted separately. The UPS is effective both in theory and in
computation.

We measure the performance of a procedure by the Hamming distance,
and use an asymptotic framework where $p \goto\infty$ and other
quantities (e.g.,~$n$, sparsity level and strength of signals) are
linked to~$p$ by fixed parameters. We find that in many cases, the UPS
achieves the optimal rate of convergence. Also, for many different~$\Omega$,
there is a common three-phase diagram in the
two-dimensional phase space quantifying the signal sparsity and signal
strength. In the first phase, it is possible to recover all
signals. In the second phase, it is possible to recover most of the
signals, but not all of them. In the third phase, successful variable
selection is impossible. UPS partitions the phase space in the same
way that the optimal procedures do, and recovers most of the signals
as long as successful variable selection is possible.

The lasso and the subset selection are well-known approaches to
variable selection. However, somewhat surprisingly, there are regions
in the phase space where neither of them is rate optimal, even in
very simple settings, such as $\Omega$ is tridiagonal, and when the
tuning parameter is ideally set.
\end{abstract}

%

\end{frontmatter}

\section{Introduction} \label{secIntro}

Consider the following sequence of regression problems:
%
%
\begin{equation} \label{model1}
Y^{(p)} = X^{(p)} \beta^{(p)} + z^{(p)},\qquad z^{(p)} \sim N(0,
I_n),\qquad n = n_p.
\end{equation}
Here, $X^{(p)}$ is an $n_p\times p$ matrix, where both $p$ and $n_p$
are large, but $p > n_p$. The $p \times1$ vector $\beta^{(p)}$ is
unknown to us, but is sparse in the sense that it has $s_p$ nonzeros
where $s_p \ll p$.
We are interested in variable selection: determining which components
of $\beta^{(p)}$ are nonzero. For notational simplicity,
we suppress the superscript $^{(p)}$ and subscript $p$ whenever there
is no confusion.

A well-known approach to variable selection is \textit{subset selection},
also known as
the $L^0$-penalization method (e.g., AIC \cite{AIC}, BIC \cite{BIC}
and RIC \cite{RIC}). This approach selects variables by minimizing the
following functional:
%
%
\begin{equation} \label{Definesubset}
\frac{1}{2} \| Y - X \beta\|_2^2 + \frac{(\lambda^{\mathrm{ss}})^2}{2} \|
\beta\|_0,
\end{equation}
where $\lambda^{\mathrm{ss}} > 0$ is a tuning parameter, and $\|\cdot\|_q$
denotes the $L^q$-norm.
The approach has good properties, but the optimization problem (\ref
{Definesubset}) is known to be
NP hard, which prohibits the use of the approach when $p$ is large.

In the middle 1990s, Tibshirani \cite{Tibshirani} and Chen et al.
\cite{Chen} proposed a trail-breaking approach which is now known
as the lasso or the basis pursuit. This approach selects variables by
minimizing a similar functional, but $\|\beta\|_0$ is replaced by
$\|\beta\|_1$.
%
%
\begin{equation} \label{Definelasso}
\tfrac{1}{2} \| Y - X \beta\|_2^2 + \lambda^{\mathrm{lasso}} \| \beta\|_1.
\end{equation}
A major advantage of the lasso is that (\ref{Definelasso}) can be
efficiently solved by the interior point method \cite{Chen}, even when
$p$ is relatively large. Additionally, in a~series of papers (e.g.,
\cite{Donoho06,DonohoTanner}), it was shown that in the noiseless
case (i.e., $z = 0$), the lasso solution is also the subset selection
solution, provided that~$\beta$ is sufficiently sparse. For these
reasons, the lasso procedure is passionately embraced by
statisticians, engineers, biologists and many others.

With that being said, an obvious shortcoming of these methods is that
the penalization term does not reflect the correlation structure in
$X$, which prohibits the method from fully capturing the essence of
the data (e.g., Zou~\cite{Zou}). However, this shortcoming is largely
due to that these methods are \textit{one-stage} procedures. This calls
for a \textit{two-stage} or \textit{multi-stage} procedure.

\subsection{Screen and clean} \label{subsecUPS}
An idea introduced in the 1960s, screen and clean, has seen a
revival recently \cite{Wasserman,FanLv}. This is a two-stage
method, where, at the first stage, we remove as many irrelevant
variables as possible while keeping all relevant ones. At the second
stage, we reinvestigate the surviving variables in hope of removing
all false positives. The screening stage has the following advantages,
some of which are elaborated in the literature:
\begin{itemize}
\item\textit{Dimension reduction}. We remove many irrelevant variables,
reducing the dimension from $p$ to a much smaller number
\cite{FanLv,Wasserman}.
\item\textit{Correlation complexity reduction}. A variable
may be correlated to many other variables, but few of which will
survive the screening; it is only correlated with a few other
surviving variables.
\item\textit{Computation complexity reduction}.
Under some conditions (e.g., Section~\ref{secUPS}), surviving
variables can be grouped into many small units, each has a~size $\leq
K$, and correlation between units is weak. These units can be
fitted separately, with computational $\mbox{cost} \leq \# \mbox{ of
units} \times2^{K}$.
\end{itemize}
Despite the perceptive vision and philosophical importance in these
works \cite{FanLv,Wasserman}, substantial vagueness remains: How to
screen? How to clean? Is screen and clean really better than the
lasso and the subset selection? This is where the Univariate
Penalization Screening (UPS) comes in.

\subsection{UPS}
The UPS is a two-stage method which contains an $U$-step and a $P$-step.
In the $U$-step, we screen with univariate thresholding
\cite{Donoho06} (also known as
marginal regression \cite{Genovese} and sure screening \cite{FanLv}).
Fix a threshold $t > 0$, and let $x_j$ be the $j$th column of $X$.
We remove the $j$th variable from the regression model if and only if
$|(x_j, Y)| < t$. The set of surviving indices is then $ {\mathcal
U}_p(t) = {\mathcal U}_p(t; Y, X) = \{j\dvtx|(x_j, Y)| \geq t, 1
\leq j \leq p \}. $

Despite its simplicity,
the $U$-step can be effective in many situations.
The key insight is that ${\mathcal U}_p(t)$ has the following important
properties:
\begin{itemize}
\item\textit{Sure Screening} (\textit{SS}). With overwhelming probability,
${\mathcal U}_p(t)$ includes all but a negligible proportion of the
signals (i.e., nonzero coordinates of $\beta$). The terminology is
slightly different from that in \cite{FanLv}.
\item\textit{Separable
After Screening} (\textit{SAS}). Define a graph where $\{1, 2, \ldots,
p\}$ is
the set of nodes, and nodes $j$ and $k$ are connected if and only if
$|(x_j, x_k)|$ is large (i.e., columns $j$ and $k$ are
``significantly'' correlated). The SAS property refers to as that
with overwhelming probability, ${\mathcal U}_p(t)$ splits into
many disconnected small-size components [a component is a maximal
connected subgraph of ${\mathcal U}_p(t)$].
\end{itemize}

We now explain how these properties pave the way for the $P$-step.
Let $\call= \{i_1, \ldots, i_K\}$ and ${\mathcal J}_0 = \{j_1, \ldots
, j_L\}$ be two subsets of $\{1,2, \ldots, p\}$, $1 \leq K$, $L \leq p$.
We have the following definition.
%
%
\begin{definition}\label{def1.1}
For any $p\times1$ vector $Y$, $Y^{\call}$ denotes the $K\times1$
vector such that $Y^{\call}(k) = Y_{i_k}$, $1 \leq k \leq K$. For any
$p\times p$ matrix $\Omega$, $\Omega^{\call, {\mathcal J}_0}$
denotes the $K\times L$ matrix such that $\Omega^{\call, {\mathcal
J}_0}(k, \ell) = \Omega(i_k, j_{\ell})$, $1 \leq k \leq K, 1 \leq
\ell\leq L$.
\end{definition}

Note that the regression model is closely related to the model
$
X' Y = X'X \beta+ X'z.
$
Restricting the attention to ${\mathcal U} = {\mathcal U}_p(t)$, we have
\[
(X'Y)^{{\mathcal U}} = (X' X \beta)^{{\mathcal U}} + (X'z)^{{\mathcal
U}} = (X'X)^{{\mathcal U}, {\mathcal V}} \beta+ (X' z)^{{\mathcal U}},
\]
where ${\mathcal V} = \{1, 2, \ldots, p\}$.
Three key observations are the following: (a)
since $z \sim N(0, I_n)$, $(X' z)^{{\mathcal U}} \sim N(0,
(X'X)^{{\mathcal U}, {\mathcal U}})$, (b) by the sure screening\vadjust{\goodbreak}
property, $(X'X)^{{\mathcal U}, {\mathcal V}} \beta\approx(X'
X)^{{\mathcal U}, {\mathcal U}} \beta^{{\mathcal U}}$ and (c) by the
SAS property, $(X'X)^{{\mathcal U}, {\mathcal U}}$ approximately equals
a block diagonal matrix, where each block corresponds to a~maximal
connected subgraph contained in ${\mathcal U}_p(t)$.
As a result, the original regression problem reduces to many small-size
regression problems that
can be solved separately, each at a modest computational cost.

In detail, fix two parameters $\lambda^{\mathrm{ups}}$ and $u^{\mathrm{ups}}$. Let
$\call= \{i_1, i_2, \ldots, i_K\} \subset{\mathcal U}_p(t)$ be a component,
and let $\mu$ be a $K\times1$ vector the coordinates of which are
either $0$ or~$u^{\mathrm{ups}}$.
Write $A=(X' X)^{\call, \call}$ for short. Let $\hat{\mu}(\call) =
\hat{\mu}(\call; Y, X, t, \lambda^{\mathrm{ups}}, u^{\mathrm{ups}}, p)$ be the
minimizer of the functional
%
%
\begin{equation} \label{PMLE}
\tfrac{1}{2} \bigl((X' Y)^{\call} - A \mu\bigr)' A^{-1} \bigl((X'
Y)^{\call} - A \mu\bigr) + \tfrac{1}{2}(\lambda^{\mathrm{ups}})^2 \| \mu\|_0.
\end{equation}
Combining all such estimates across different components of ${\mathcal
U}_p(t)$ gives the UPS estimator, denoted by $\hb^{\mathrm{ups}} = \hb
^{\mathrm{ups}}(Y, X; t, \lambda^{\mathrm{ups}}, u^{\mathrm{ups}}, p)$,
\[
\hb^{\mathrm{ups}}_j = \cases{
(\hat{\mu}(\call))_k, &\quad if $j = i_k \in\call$ for some $\call= \{
i_1, i_2, \ldots, i_K\} \subset{\mathcal
U}_p(t)$,\cr
0, &\quad if $j \notin{\mathcal U}_p(t_p)$.}
\]

The UPS uses three tuning parameters $(t, \lambda^{\mathrm{ups}}, u^{\mathrm{ups}})$.
In many cases, the performance
of the UPS is relatively insensitive to the choice of $t$, as long as
it falls in a certain range. The parameter $\lambda^{\mathrm{ups}}$
has a similar role to those of the lasso and the subset selection, but
there is a major difference:
the former can be conveniently estimated using the data, whereas how to
set the latter remains an open problem. See Section \ref{secUPS} for
more discussion.

We are now ready to answer the
questions raised in the end of Section~\ref{subsecUPS}: UPS indeed has
advantages over the lasso and the subset selection.
In Sections \ref{subsecBounds}--\ref{subseclasso},
we establish a theoretic framework and investigate these procedures closely.
The main finding is the following: for a wide range of design matrices
$X$, the Hamming distance of the UPS achieves the optimal rate of
convergence. In contrast, the lasso and the subset selection may be
rate nonoptimal, even for very simple design matrices.

\subsection{Sparse signal model and universal lower bound}
\label{subsecBounds}

We model $\beta$ by
%
%
\begin{equation} \label{betaadd}
\beta_j \stackrel{\mathrm{i.i.d.}}{\sim} (1 - \eps) \nu_0 +
\eps\pi,\qquad
0 < \eps< 1, 1 \leq j \leq p,
\end{equation}
where $\nu_0$ is the point mass at $0$, and $\pi$ is a distribution
that has no mass at~$0$. We
use $p$ as the driving asymptotic parameter and
allow $(\eps, \pi)$ to depend on $p$.
Fix $0 < \vartheta< 1$ and recall that $s_p$ is the number of signals.
We calibrate
%
%
\begin{equation} \label{Defineeps}
\eps= \eps_p = p^{-\vartheta}\qquad \mbox{so that $s_p \sim p
\eps_p = p^{1 - \vartheta}$}.
\end{equation}
For any variable selection procedure $\hat{\beta} = \hb(Y | X)$, we measure
the loss by the Hamming distance
\[
h_p(\hb, \beta| X) = h_p(\hb, \beta; \eps_p, \pi_p, n_p | X) =
E_{\eps_p, \pi_p} \Biggl[ \sum_{j =1}^p 1\bigl(\sgn(\hb_j) \neq
\sgn(\beta_j) \bigr) \Biggr],\vadjust{\goodbreak}
\]
where $\sgn(0)=0$. In the context of variable selection, the Hamming
distance is a natural choice for loss function. While the focus of this
paper is on selection error where we use $L_0$-loss, the idea can be
extended to the estimation setting where we use $L_q$-loss ($0 < q <
\infty)$, but we have to perform an additional step of least square
fitting after the selection.

Somewhat surprisingly, there is a lower bound for the Hamming distance
that holds for all sample size $n$ and design matrix $X$ (and so
``universal lower bound''). The following notation is frequently used
in this paper.
%
%
\begin{definition}\label{def1.2}
$L_p > 0$ is a multi-$\log(p)$ term which may change from occurrence
to occurrence, such that for any fixed $\delta> 0$, $\lim
_{p\rightarrow\infty} L_p \cdot p^\delta=\infty$ and $\lim
_{p\rightarrow\infty} L_p p^{-\delta}=0$.
\end{definition}

Now, fixing $r > 0$, we introduce
%
%
\begin{equation} \label{Definetau}
\tau_p = \tau_p(r) = \sqrt{2 r \log p}
\end{equation}
and
$\lambda_p = \lambda_p(\eps_p, \tau_p) = \frac{1}{\tau_p} [
\log(\frac{1 - \eps_p}{\eps_p}) + \frac{\tau_p^2}{2} ]$.
Let $\bar{\Phi} = 1 - \Phi$ be the survival function of $N(0,1)$.
The following theorem is proved in
\cite{UPSSupp}.
%
%
\begin{theorem}[(Lower bound)] \label{thmLB}
Fix $\vartheta\in(0,1)$, $r > 0$ and a
sufficiently large~$p$. Let $\eps_p$, $s_p$ and $\tau_p$ be as in
(\ref{Defineeps}) and (\ref{Definetau}), and
suppose the support of~$\pi_p$ is contained in $[-\tau_p, 0) \cup(0,
\tau_p]$. For any
fixed $n$ and matrix \mbox{$X = X^{(p)}$} such that $X'X$ has unit diagonals,
$h_p(\hat{\beta}, \beta| X) \geq
s_p \cdot[(1-\varepsilon_p) \bar{\Phi}(\lambda_p)/\varepsilon_p+\break\Phi
(\tau_p-\lambda_p)]$.
\end{theorem}

Note that as $p \goto\infty$,
%
%
\begin{equation} \label{rate}
\frac{1-\varepsilon_p}{\varepsilon_p} \bar{\Phi}(\lambda_p) + \Phi
(\tau_p-\lambda_p) \geq \cases{ L_p \cdot p^{- (r - \vartheta)^2/(4r)},
&\quad $r > \vartheta$, \vspace*{2pt}\cr \bigl(1 + o(1)\bigr), &\quad
$r < \vartheta$.}
\end{equation}
It may seem counterintuitive that the lower bound does not depend
on~$n$, but this is due to the way we normalize $X$. In the case of
orthogonal design [i.e., coordinates of $X$ and i.i.d. from $N(0,
1/n)$], the lower bound can be achieved
by either the lasso or marginal regression \cite{Genovese}. Therefore,
the orthogonal design is among the best in terms of
the error rate.

Theorem \ref{thmLB} says that if we have $p^{1 - \vartheta}$ signals,
and the maximal signal strength is slightly smaller than $\sqrt{2
\vartheta\log(p})$, then the Hamming distance of any procedure
cannot be substantially smaller than $s_p$, and so successful variable
selection is impossible. In the sections below, we focus on the case
where the signal strength is larger than $\sqrt{2 \vartheta\log
(p)}$, so that successful variable selection is possible.

The universality of the lower bound hints it may not be tight for
nonorthogonal~$X$.
Fortunately, it turns out that in many interesting cases, the lower
bound is tight.
To facilitate the analysis, we invoke the random design model.

\subsection{Random design, connection to Stein's normal means model}
Write
$
X = (x_1, x_2, \ldots, x_p) = (X_1, X_2,\ldots, X_n)'$.
We model $X_i$ as i.i.d. samples from a $p$-variate zero-mean Gaussian
distribution,
%
%
\begin{equation} \label{modelX}
X_i \stackrel{\mathrm{i.i.d.}}{\sim} N\biggl(0, \frac{1}{n} \Omega\biggr).
\end{equation}
The $p\times p$ matrix $\Omega=\Omega^{(p)}$ is unknown, but for
simplicity we assume it has unit diagonals. The normalizing constant
$1/n$ is chosen so that the diagonals of the Gram matrix $X'X$ are
approximately $1$.
Fixing $\theta\in(1 - \vartheta, 1)$, we let
%
%
\begin{equation} \label{Definenp}
n = n_p = p^{\theta}.
\end{equation}
Note that $s_p \ll n_p \ll p$ as $p \goto\infty$. For successful
variable selection, it is almost necessary to have $s_p \ll n_p$
\cite{Donoho06}. Also, denoting the distribution of~$X$ by $F = F_p$,
note that for any variable selection procedure, the \textit{overall
Hamming distance} is $\hamm_p(\hb, \beta) = E_F[ h_p(\hb|X)]$.

Model (\ref{modelX}) is called the \textit{random design model} which may
be found in the following application areas:
\begin{itemize}
\item\textit{Compressive sensing.} We are interested in a
$p$-dimensional sparse vector $\beta$. We measure $n$ general linear
combinations of $\beta$ and then reconstruct it. For $1 \leq i \leq
n$, choose a $p\times1$ coefficient vector $X_i$, and observe $Y_i
= X_i' \beta+ z_i$, where $z_i \sim N(0, \sigma^2)$ is noise. For
computational and storage concerns, one usually chooses $X_i$'s as
simple as possible. Popular choices of $X_i$ include Gaussian design,
Bernoulli design, circulant design, etc. \cite{Donoho06,NowakCS}.
Model (\ref{modelX}) belongs to Gaussian design.
\item\textit{Privacy-preserving data mining.} The vector $\beta$ may
contain some
confidential information (e.g., HIV-diagnosis results of a community)
that we must protect. While we cannot release the whole vector, we
must allow data mining to some extent, because, for example, the study
is of public interest and is supported by federal funding. To
compromise, we allow queries as follows. For each query, the database
randomly generates a $p\times1$ vector $X_i$, and releases both $X_i$
and $Y_i = X_i' \beta+ z_i$ to the querier, where $z_i \sim N(0,
\sigma^2)$ is a noise term. For privacy concerns, the number of
allowed queries is much smaller than $p$. Popular choices of $X_i$
include Gaussian design and Bernoulli design \cite{Nissim}.
\end{itemize}

Random design model is closely related to Stein's normal means model
$W \sim N(\beta, \Sigma)$,
where $\Sigma= \Omega^{-1}$.
To see the point,
recall that model (\ref{model1}) is closely related to the model $X' Y
= X' X \beta+ X' z$.
Since the rows of $X$ are i.i.d. samples from $N(0, \frac{1}{n} \Omega
)$ and $s_p \ll n_p \ll p$, we expect to see that
$X' X \beta\approx\Omega\beta$ and $X' z \approx N(0, \Omega)$,
and so that $X' Y \approx N(\Omega\beta, \Omega)$.
Therefore, Stein's normal means model can be viewed as an idealized
version of the random design model.
This suggests that solving the variable selection problem opens doors
for solving Stein's normal means problem, and vice versa.

\subsection{Optimality of the UPS}
The main results of this paper are Theorems \ref{thmUB} and \ref
{thmadaptive} in Section \ref{secUPS}.
To state such results, we need relatively long
preparations. Therefore, we sketch these results below, but leave the
formal statements to later.
In models (\ref{model1}), (\ref{betaadd}) and (\ref{modelX}), let
$(s_p, \tau_p, n_p)$ be as in~(\ref{Defineeps}), (\ref{Definetau})
and (\ref{Definenp}). Suppose:
\begin{itemize}
\item Each row of $\Omega$ satisfies a certain summability condition,
so it has relatively few large coordinates.
\item The support of
$\pi_p$ is contained in $[\tau_p, (1 + \eta) \tau_p]$, where $\tau
_p =
\sqrt{2 r \log(p)}$, and $\eta$ is a constant to be defined later. We
suppose $r > \vartheta$, so that successful variable selection is
possible; see Theorem \ref{thmLB}.
\item Either all coordinates of
$\Omega$ are positive, or that $r/\vartheta\leq3 + 2 \sqrt{2}$ (so
that we won't have too many ``signal cancellations'' \cite{Wasserman}).
\end{itemize}
Fix $0 < q \leq(\vartheta+ r)^2/(4r)$, and set the tuning parameters
$(t, \lambda^{\mathrm{ups}}, u^{\mathrm{ups}})$ by
\[
t_p^* = t_p^*(q) = \sqrt{2 q \log p},\qquad \lambda^{\mathrm{ups}} = \lambda
_p^{\mathrm{ups}} = \sqrt{2 \vartheta\log(p)},\qquad
u^{\mathrm{ups}} = u_p^{\mathrm{ups}} =
\tau_p.
\]
The main result is that, as $p \goto\infty$, the ratio between the
Hamming error of the UPS and $s_p$ is no grater than $ L_p
p^{-(\vartheta- r)^2/(4r)}$. Comparing this with Theorem \ref{thmLB}
gives that the lower bound is tight, and the UPS
is rate optimal.

\subsection{Phase diagram for high-dimensional variable selection}
\label{subsecPhase}
The above results reveal a watershed phenomenon as follows. Suppose we
have roughly $s_p = p^{1-\vartheta}$ signals. If
the maximal signal strength is slightly smaller than $\sqrt{2
\vartheta\log p}$, then the Hamming distance of
any procedure cannot be substantially smaller than~$s_p$, hence
successful variable selection is impossible. If the minimal signal
strength is slightly larger than $\sqrt{2 \vartheta\log p}$,
then there exist procedures (UPS is one of them) whose Hamming
distances are substantially smaller than $s_p$, and they manage to
recover most signals.

\begin{figure}

\includegraphics{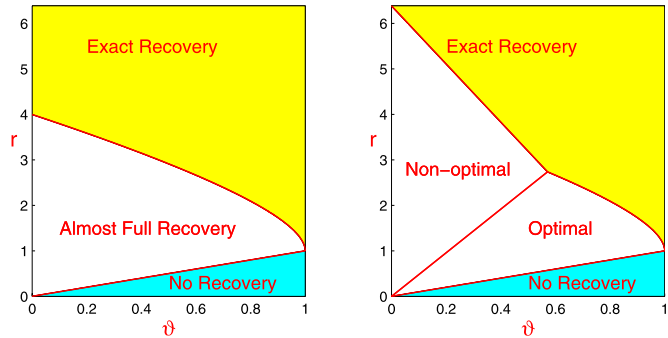}

\caption{Left: phase diagram.
In the yellow region, the UPS recovers all signals with high
probability. In the white region, it is possible (i.e., UPS) to recover
almost all signals, but impossible to recover all of them. In the cyan
region, successful variable selection is impossible.
Right: partition of the phase space by the lasso for the tridiagonal
model~(\protect\ref{tridiagonal1})--(\protect\ref{tridiagonal2})
($a =
0.4$). The lasso is rate nonoptimal in the nonoptimal region. The
region of exact recovery by the lasso is substantially smaller than
that displayed on the left.} \label{figLassoPhase}
\end{figure}

The phenomenon is best described in the special case where
$\pi_p = \nu_{\tau_p}$ is the point mass at $\tau_p$,
with $\tau_p = \sqrt{2 r \log p}$ as in (\ref{Definetau}).
If we call the two-dimensional domain $\{(\vartheta, r)\dvtx0<
\vartheta< 1, r > 0\}$ the \textit{phase space}, then
the theorems say that the phase space is partitioned into three regions:
\begin{itemize}
\item\textit{Region of no recovery}
($0 < \vartheta< 1$, $0 < r <
\vartheta$). In this region, the Hamming distance of any $\mbox{procedure}
\gtrsim s_p$, and successful variable selection is impossible.
\item\textit{Region of almost full recovery} [$0 < \vartheta< 1$,
$\vartheta< r < (1 + \sqrt{1 - \vartheta})^2$]. In this region,
there are procedures (e.g., UPS) whose Hamming errors are much larger
than $1$, but are also much smaller than $s_p$. In this region, it is
possible to recover most of the signals, but not all of them.
\item\textit{Region of exact recovery} [$0 < \vartheta< 1$, $r > (1 +
\sqrt{1 - \vartheta})^2$]. In this region, there are procedures
(e.g., UPS) that recover all signals with probability~\mbox{$\approx1$}.
\end{itemize}
See Figure \ref{figLassoPhase} (left panel) for these regions.
Note that the partitions are the same for many choices of $\Omega$.
Because of the partition of the phases, we call this the phase diagram.
The UPS is optimal in the sense that it partitions the phase space in
exactly the same way as do the optimal procedures.

%

The phase diagram provides a benchmark for variable selection.
The lasso would be optimal if it partitions the phase space
in the same way as in the left panel of Figure \ref{figLassoPhase}.
Unfortunately, this is not the case, even for very simple~$\Omega$.
Below we investigate the case where $X'X$ is a tridiagonal matrix, and
identify precisely the regions where the lasso is rate optimal and
where it is rate nonoptimal.
More surprisingly, there is a region in the phase space where the
subset selection is also rate nonoptimal.

\subsection{Nonoptimal region for the lasso} \label{subseclasso}
In Sections \ref{subseclasso} and \ref{subsecsubset}, we temporarily
leave the random design model and consider Stein's normal means model,
which is an idealized version of the former. Using an idealized version
is mainly
for mathematical convenience, but the gained insight is valid
in much broader settings:
if a procedure is nonoptimal in simple cases,
we should not expect them to be optimal in more complicated cases.

In this spirit, we consider Stein's normal means model
%
%
\begin{equation} \label{tridiagonal1}
\tilde{Y} \equiv X' Y \sim N(\Omega\beta, \Omega),
\end{equation}
where $\beta$ is as in (\ref{betaadd}) with $\tau_p = \nu_{\pi_p}$
and $\pi_p = \sqrt{2 r \log(p)}$.
To further simplify the study, we fix $a \in(0, 1/2)$ and take $\Omega
$ as the tridiagonal matrix $T(a)$:
%
%
\begin{equation} \label{tridiagonal2}
T(a)(i,j) = 1\{i = j\} + a \cdot1\{ |i -j | = 1\},\qquad
1 \leq i, j \leq p.
\end{equation}
Note that in this case the UPS partitions the phase space optimally.

We now discuss the phase diagram of the lasso.
The region $\{(\vartheta, r)\dvtx0 < \vartheta< 1, r > \vartheta\}$
is partitioned into three regions as follows (see Figure \ref{figLassoPhase}):
\begin{itemize}
\item\textit{Nonoptimal region}: $0 < \vartheta< 2a(1 + a)^{-1}$ and
$\frac{1}{a} (1 + \sqrt{1 - a^2}) \vartheta< r < (1 + \sqrt
{\frac{1 + a}{1 - a}} )^2 (1 - \vartheta)$.
In this region, the lasso is rate nonoptimal [i.e., the Hamming
distance is $L_p \cdot p^{c}$ with constant $c > 1 - (\vartheta+
r)^2/(4r)$], even when the tuning parameter is set ideally.
\item\textit{Optimal region}: $0 < \vartheta< 1$ and $\vartheta< r <
\frac{1}{a}(1 + \sqrt{1 - a^2})\vartheta$ and $r < (1 + \sqrt{1 -
\vartheta})^2$. In this region, if additionally $a \geq1/3$, then the
lasso may be rate optimal if the tuning parameter is set ideally. The
discussion on the case $0 < a < 1/3$ is
tedious so we skip it.
\item\textit{Region of exact recovery}: $0 < \vartheta< 1$ and $r > (1 +
\sqrt{1 - \vartheta})^2$ and $r > (1 + \sqrt{\frac{1 + a}{1 -
a}} )^2 (1 - \vartheta)$. In this region, if the tuning
parameter is set ideally, the lasso may
yield exact recovery with high probability. Region of exactly recovery
by the lasso is substantially smaller
than that of the UPS. There is a sub-region in the phase space where
the UPS yields exact recovery, but the lasso could not even when the
tuning parameter is set ideally.
\end{itemize}

For discussions in the case where $\Omega$ is the identity matrix,
compare \cite{Genovese,Wainwright}. The above results are proved in
Theorem \ref{thmlasso}, where we derive a lower bound for the Hamming
errors by the lasso. In \cite{Thesis}, we show that the lower bound
is tight for properly large $\vartheta$, but is not when $\vartheta$
is small. It is, however, tight for all $\vartheta\in(0,1)$ if we
replace model (\ref{betaadd}) by a closely related model, namely
(2.2) and (2.3) in~\cite{HallJin}. For these reasons, the nonoptimal
region of the lasso may be larger than that illustrated in Figure
\ref{figLassoPhase}. The discussion on the exact optimal rate of
convergence for the lasso is tedious and we skip it.

Why is the lasso nonoptimal? To gain insight, we introduce the term of
\textit{fake signal}, a noise coordinate that may look like a signal due
to correlation.

%
\begin{definition}\label{def1.3}
We say\vspace*{1pt} that $\tilde{Y}_j$ is a signal if $\beta_j \neq0$, is a fake
signal if $(\Omega\beta)_j \neq0$ and $\beta_j = 0$, and is a
(pure) noise if
$\beta_j = (\Omega\beta)_j = 0$.
\end{definition}

With the tuning parameter set ideally, the lasso is able to distinguish
signals from pure noise, but it does not filter out fake signals efficiently.
In the optimal region of the lasso, the number of falsely kept fake
signals is much smaller than the optimal rate, so it is negligible; in
the nonoptimal region, the number becomes much larger than the optimal
rate, and so is nonnegligible.
This suggests that when $X'X$ moves away from the tridiagonal case, the
partitions of the regions by the lasso may change, but the nonoptimal
region of the lasso continues to exist in rather general situations.

The nonoptimality of the lasso is largely due to the fact that it is a
one-stage method.
An interesting question is whether UPS continues to work well if we
replace the univariate thresholding by the lasso in the screening stage.\vadjust{\goodbreak}
The disadvantage of this proposal is that, compared to the univariate
thresholding, the lasso is both slower in computation and harder to
analyze in theory.
Still, one would hope the lasso could perform well in screening.

With that being said, we note that the implementation of the lasso only
needs minimal assumption on the model, which makes it very attractive,
especially in complicated situations. In comparison, we need both
signal sparsity and graph sparsity to implement the UPS, and how to
extend it to more general settings remains unknown. The exploration
along this line is continued in our forthcoming manuscripts
\cite{JinZhang1,JinZhang2,FJZ}; see details therein.

\subsection{Nonoptimal region for the subset selection} \label{subsecsubset}
\begin{figure}[b]

\includegraphics{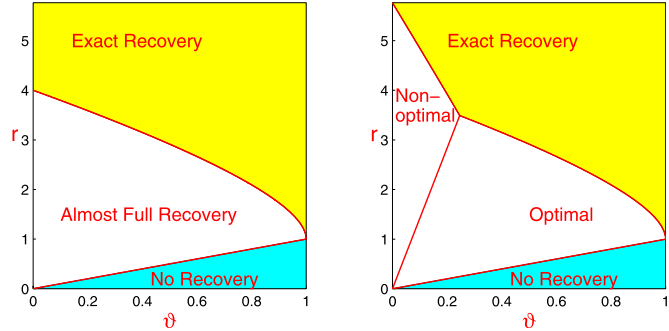}

\caption{Left: a re-display of the left panel of Figure
\protect\ref{figLassoPhase}. Right: partition of the phase space by
the subset selection in the tridiagonal model
(\protect\ref{tridiagonal1})--(\protect\ref{tridiagonal2}) ($a = 0.4$).
The subset selection is not rate optimal in the nonoptimal region. The
exact recovery region by the subset selection is substantially smaller
than that of the optimal procedure, displayed on the left. }
\label{figSSPhase}
\end{figure}

The discussion on the subset selection is similar to that for the lasso
so we keep it brief.
Introduce
$v_1(a) = \frac{2 - \sqrt{1 - a^2}}{\sqrt{1 - a^2} (1 - \sqrt{1 - a^2})}$
and $ v_2(a) = 2 \sqrt{1 - a^2} -1$.
Similarly, the phase space partitions into three regions as follows:
\begin{itemize}
\item\textit{Nonoptimal region}: $0 < \vartheta<\frac{4
v_1(a)}{(v_1(a)+1)^2}$ and $v_1(a) \vartheta< r < [\frac
{1}{v_2(a)} (\sqrt{1 - 2 \vartheta} + \sqrt{1 - 2 \vartheta+
\vartheta v_2(a)} ) ]^2$.
\item\textit{Optimal region}: $0 < \vartheta< 1$ and
$\vartheta< r < v_1(a) \vartheta$ and $r < (1 + \sqrt{1 - \vartheta})^2$.
\item\textit{Exact recovery region}: $0 < \vartheta< 1$, $r > (1 + \sqrt
{1 - \vartheta})^2$ and $r > [\frac{1}{v_2(a)} (\sqrt{1
- 2 \vartheta} + \sqrt{1 - 2 \vartheta+ \vartheta v_2(a)}
)]^2$.
\end{itemize}
See Theorem \ref{thmsubset} for proofs and Figure \ref{figSSPhase}
for illustration. Similar to the remarks in Section \ref{subseclasso}, the
region of exact recovery and the optimal region of the subset selection
may be smaller than those illustrated in Figure \ref{figSSPhase}.

%

The reason why the subset selection is nonoptimal is almost the
\textit{opposite} to that of the lasso: the lasso is nonoptimal for it
is too
loose on fake signals, but the subset selection is nonoptimal for it
is too harsh on signal clusters (pairs/triplets, etc.). With the tuning
parameter set ideally, the subset selection is effective in filtering
out fake signals, but it also tends to kill one or more signals when
the true signals appear in clusters. These falsely killed signals
account for the nonoptimality.
See Section \ref{secSS} for details.

\subsection{Connection to recent literature} \label{subsecconnection}
This work is related to recent literature on oracle property
\cite{Zou,Meinshausen}, but is different in important ways. A~procedure
has the oracle property if it yields exact recovery. However, exact
recovery is rarely seen in applications, especially when $p \gg n$. In
many applications (e.g., genomics), a~large $p$ usually means that
signals are sparse or rare, and a small $n$ usually means signals are
weak. For rare and weak signals, exact recovery is usually impossible.
Therefore, it is both scientifically more relevant and technically more
challenging to compare error rates of different procedures than to
investigate when they satisfy the oracle property.

The work is also related to \cite{Candes,Zhang} on asymptotic
minimaxity, where the lasso was shown to be asymptotic rate optimal in
the worst-case scenario. While their results seem to contradict with
those in this paper, the difference can be easily reconciled. In the
minimax approach, the asymptotic least favorable distribution of
$\beta$ is given by $\beta_j \stackrel{\mathrm{i.i.d.}}{\sim} (1 - \eps_p)
\nu_0 + \eps_p \nu_{\tau_p}$, where \mbox{$\eps_p = p^{-\vartheta}$},
$\tau_p = \sqrt{2 r \log p}$ and notably $\vartheta= r$, which
corresponds the boundary line of the region of no recovery in the
phase space (e.g.,~\cite{Zhang}, pages~18 and~19,~\cite{ABDJ},
Section 3).
This suggests that the minimax approach has limitations: it reduces the
analysis to the worst-case scenario, but the worst-case scenario may be
outside the range of interest.
In our approach, we let $(\vartheta, r)$ range freely, and evaluate a
procedure based on how it partitions the phase space. Our approach has
a similar spirit to that in \cite{DonohoTanner}.

The work is also related to the adaptive lasso \cite{Zou}. The
adaptive lasso is similar to the lasso, but the $L^1$-penalty $\lambda
^{\mathrm{lasso}} \|\beta\|_1$ is replaced by the weighted
$L^1$-penalty $\sum_{j = 1}^p w_j |\beta_j|$, where $w = (w_1, \ldots
, w_p)'$ is the weight vector.
Philosophically, we can view the adaptive lasso as a screen and clean
method. Still, the proposed approach is
different from
the adaptive lasso in important ways. First, Zou \cite{Zou} suggested
weight choices by the least squares estimate, which is only feasible
when $p$ is small.
In fact, when $p \gg n$, our results suggest that feasible weights
should be very sparse, while the weights suggested by the least squares
estimates are usually dense.
Second, for the surviving indices, we first partition them into many
disjoint units of small sizes, and then fit them individually. The
adaptive lasso fits all surviving variables together, which is
computationally more expensive. Last, we use penalized MLE in the clean
step while the adaptive lasso uses $L^1$-penalty. As pointed out
before, the $L^1$-penalty in the clean step is too loose on fake
signals, which prohibits the procedure from being rate optimal.

The work is also related to other multi-stage methods, for example,
the threshold lasso \cite{zhouthresholded2010} or the LOL
\cite{kerkyacharianlearning2009}. These methods first use the lasso
and the OLS for variable selection, respectively, followed by an
additional thresholding step. However, by an argument similar to that
in Sections~\ref{subseclasso} and~\ref{subsecsubset}, it is not hard to
see that these procedures do not partition the phase diagram optimally.

\subsection{Contents} \label{subseccontent}
In summary, we propose the UPS as a two-stage method for variable selection.
We use Univariate thresholding in the screening step for its
exceptional convenience in computation, and we use penalized MLE in the
cleaning step because it is the only procedure we know so far that
yields the optimal rate of convergence.
On the other hand, the lasso and even the subset selection do not
partition the phase space optimally.

The remaining sections are organized as follows. Section \ref{secUPS}
discusses the UPS procedure and
the upper bound for the rate of
convergence. The section also addresses how to estimate the tuning
parameters of the UPS and the convergence rate of the resultant
plug-in procedure. Section \ref{secrefinement} discusses a
refinement of the UPS for moderately large $p$. Section \ref{seclasso}
discusses the behavior of the lasso and the subset selection. Section
\ref{secSimul} discusses numerical results where we compare the UPS
with the lasso (the subset selection is computationally infeasible for
large $p$ so is not included for comparison). Due to limited space, we
do not include proofs in this paper. The proofs can be found in the
supplementary material for the paper \cite{UPSSupp}.

Below is some notation we use in this paper.
Fix $0 < q < \infty$. For a $p\times1$ vector~$x$, $\|x\|_q$ denotes
the $L^q$-norm of $x$, and we omit the subscript when $q = 2$. For a
$p\times p$ matrix $M$, $\|M\|_q$ denotes the matrix $L^q$-norm, and $\|
M\|$ denotes the spectral norm.

\section{UPS and upper bound for the Hamming distance} \label{secUPS}

In this section, we establish the upper bound for the Hamming distance
and show that the UPS is rate optimal.
We begin by discussing necessary notation. We then discuss the
$U$-step and its sure screening and SAS properties.
Next, we show how the regression problem reduces to many separate
small-size regression problems and explain the rationale of using the
penalized MLE in the $P$-step.
We conclude the section by the rate optimality of the UPS, where the
tuning parameters are either set ideally or estimated.

Since different parts of our model are introduced separately in
different subsections, we summarize them as follows.
The model we consider is
%
%
\begin{equation} \label{modelrecap}
Y = X \beta+ z,\qquad z \sim N(0, I_n),
\end{equation}
where
%
%
\begin{eqnarray} \label{Assumptiona}
X_i &\stackrel{\mathrm{i.i.d.}}{\sim}& N\biggl(0,
\frac{1}{n}\Omega\biggr),\nonumber\\[-8pt]\\[-8pt]
\beta_j
&\stackrel{\mathrm{i.i.d.}}{\sim}& (1 - \eps_p) \nu_0 + \eps_p
\pi_p,\qquad
1 \leq i \leq n, 1 \leq j \leq p.\nonumber
\end{eqnarray}
Fixing $\theta> 0$, $\vartheta> 0$, and $r > 0$, we calibrate
%
%
\begin{equation} \label{Assumptionb}
\eps_p = p^{-\vartheta}, \qquad\tau_p = \sqrt{2 r \log p},\qquad
n_p = p^{\theta},
\end{equation}
assuming that
%
%
\begin{equation} \label{Assumptionc}
\theta< (1 - \vartheta).
\end{equation}
Recall that the optimal rate of convergence is $L_p p^{1 -(\vartheta+
r)^2/(4r)}$. In this section, we focus on the case where the exponent
$1 - (\vartheta+ r)^2/(4r)$ falls between~$0$ and $(1 - \vartheta)$,
or equivalently,
%
%
\begin{equation} \label{Assumptiond}
\vartheta< r < \bigl(1 + \sqrt{1 - \vartheta}\bigr)^2.
\end{equation}
In the phase space, this corresponds to the region of almost full
recovery. The case $r < \vartheta$ corresponds to
the region of no recovery and is studied in Theorem~\ref{thmLB}. The
case $r > (1 + \sqrt{1 - \vartheta})^2$ corresponds
to the region of exact recovery. The discussion in this case is similar
but is much easier, so we omit it.

Next, fixing $A > 0$ and $\gamma\in(0, 1)$, introduce
\[
{\mathcal M}_p(\gamma, A) = \Biggl\{ \Omega\dvtx p\times p
\mbox{ correlation matrix}, \sum_{j = 1}^p |\Omega(i,j)|^{\gamma} \leq
A, \forall 1 \leq i \leq p \Biggr\}.
\]
For any $\Omega$, let $U = U(\Omega)$ be the $p \times p$ matrix
satisfying $U(i,j) = \Omega(i,j) 1 \{i < j\}$, and let
$d(\Omega) = \max\{ \| U(\Omega)\|_1, \|U(\Omega)\|_{\infty} \}$.
Fixing $\omega_0 \in(0,1/2)$, introduce
${\mathcal M}_p^*(\omega_0, \gamma, A) = \{ \mbox{$\Omega\in
{\mathcal M}_p(\gamma, A)\dvtx d(\Omega) \leq\omega_0$} \}$,
and a subset of ${\mathcal M}_p^*(\omega_0, \gamma, A)$,
\[
{\mathcal M}_p^+(\omega_0, \gamma, A) = \{ \Omega\in{\mathcal
M}_p^*(\omega_0, \gamma, A)\dvtx\Omega(i,j) \geq0 \mbox{ for
all } 1\leq i, j \leq p \}.
\]
For any $\Omega\in{\mathcal M}_p^*(\omega_0, \gamma, A)$, the
eigenvalues are contained in $(1 - 2 \omega_0, 1 + 2 \omega_0)$, so
$\Omega$ is positive definite (when
$\omega_0 > 1/2$, $\Omega$ may not be positive definite).

Last, introduce a constant $\eta= \eta(\vartheta, r, \omega_0)$ by
%
%
\begin{equation} \label{Defineconstantc1}
\eta= \frac{\sqrt{\vartheta r}}{(\vartheta+ r)\sqrt{1 + 2 \omega
_0}} \min\biggl\{ \frac{2 \vartheta}{r}, 1 - \frac{\vartheta}{r},
\sqrt{2(1 - \omega_0)} -1 + \frac{\vartheta}{r} \biggr\}.
\end{equation}
We suppose the support of signal distribution $\pi_p$ is contained in
%
%
\begin{equation} \label{Assumptione}
[\tau_p, (1 + \eta) \tau_p],
\end{equation}
where $\tau_p = \sqrt{2 r \log(p)}$ as in (\ref{Definetau}). This
assumption is only needed for proving the main lemma of the $P$-step
(Lemma A.5, \cite{UPSSupp}) and can be relaxed for proving other
lemmas. Also, we assume the signals are one-sided mainly for
simplicity. The results can be extended to the case with two-sided
signals.

We now discuss the $U$-step. As mentioned before, the benefits of the
$U$-step are threefold: dimension reduction, correlation complexity
reduction, and computation cost reduction. The $U$-step is able to
achieve these goals simultaneously because it satisfies the sure
screening property and the SAS property, which we now discuss separately.

\subsection{The sure screening property of the $U$-step}
Recall that in the $U$-step, we remove the $j$th variable if and only
if $|(x_j, Y)| < t$ for some threshold $t > 0$. For simplicity,
we make a slight change and remove the $j$th variable if and only if
$
(x_j, Y) < t.
$
When the signals are one-sided, the change makes negligible difference.
Fixing a constant $q \in(0, (\vartheta+ r)^2/(4r))$,
we set the threshold $t$ in the $U$-step
%
%
\begin{equation} \label{Definet}
t^*_p = t_p^*(q) = \sqrt{2 q \log(p)}.
\end{equation}

%
%
\begin{lemma}[(Sure screening)] \label{lemmascreen}
In model (\ref{modelrecap})--(\ref
{Assumptiona}), suppose (\ref{Assumptionb})--(\ref{Assumptione}) hold,
and $t^*_p$ is as in (\ref{Definet}). For sufficiently large $p$, if
$\Omega^{(p)} \in{\mathcal M}_p^+(\omega_0, \gamma, A)$, then as $p
\goto\infty$, $\sum_{j = 1}^p P(x_j' Y < t_p^*, \beta_j \neq0)
\leq L_p p^{1 -{(\vartheta+ r)^2}/({4r})}$.
The claim remains true if alternatively $\Omega^{(p)} \in{\mathcal
M}_p^*(\omega_0, \gamma, A)$, but $r/\vartheta\leq3 + 2 \sqrt{2}$.
\end{lemma}

This says that the Hamming errors we make in the $U$-step are not
substantially larger than the optimal rate of convergence, and thus negligible.

\subsection{The SAS property of the $U$-step}

We need some terminology in graph theory (e.g., \cite{Diestel}).
A graph $G = (V, E)$ consists of two finite sets $V$ and $E$, where $V$
is the set of \textit{nodes}, and
$E$ is the set of \textit{edges}.
%
A \textit{component}~$\call$ of $V$ is a maximal connected subgraph,
denoted by $\call\lhd V$.
For any node $v \in V$,
there is a unique component $\call$ such that $v \in\call\lhd V$.

Fix a $p\times p$ symmetric matrix $\Omega_0$ which is presumably
sparse. If we let $V_0 = \{1, 2, \ldots, p\}$ and
say nodes $i$ and $j$ are \textit{linked} if and only if $\Omega_0(i,j)
\neq0$, then we have
a graph $G = (V_0, \Omega_0)$.
Fix $t > 0$.
Recall that ${\mathcal U}_p(t)$ is the set of surviving indices in the $U$-step
%
%
\begin{equation} \label{DefineU}
{\mathcal U}_p(t) = {\mathcal U}_p(t, Y, X) = \{j\dvtx(x_j, Y) \geq t,
1 \leq j \leq p\}.
\end{equation}
Note that the induced graph $({\mathcal U}_p(t), \Omega_0)$ splits
into many components.
%
%
\begin{definition}
Fix an integer $K \geq1$. We say that ${\mathcal U}_p(t)$ has the
separable after screening (SAS) property with respect to $(V_0, \Omega
_0, K)$ if each component of the graph $({\mathcal U}_p(t), \Omega_0)$
has no more than $K$ nodes.
\end{definition}

Note that if ${\mathcal U}_p(t)$ has the SAS property with respect to
$(V_0, \Omega_0, K)$. Then for all $s > t$, ${\mathcal U}_p(s)$ also
has the SAS property with respect to $(V_0, \Omega_0, K)$.

Return to model (\ref{modelrecap})--(\ref{Assumptiona}). We hope to
relate the regression setting to a graph $(V_0, \Omega_0)$, and use it
to spell out the SAS property. Toward this end, we set
$V_0 = \{1, 2, \ldots, p\}$.
As for $\Omega_0$, a natural choice is the matrix $\Omega$ in (\ref
{Assumptiona}). However, the SAS property makes more sense if $\Omega
_0$ is sparse and known, while $\Omega$ is neither. In light of this,
we take $\Omega_0$ to be a regularized empirical covariance
matrix.

In detail, let $\hat{\Omega} = X'X$ be the empirical covariance matrix.
Recall that
$X = (X_1, X_2, \ldots, X_n)'$ and
$ X_i \sim N(0, \frac{1}{n} \Omega)$.
It is known \cite{Bickel2008} that there\vadjust{\goodbreak} is a~constant $C > 0$ such
that with probability $1 - o(1/p^2)$, for all $1 \leq i, j \leq p$,
%
%
\begin{equation} \label{Bickel}
|\hat{\Omega}(i,j) - \Omega(i,j) | \leq C \sqrt{\log(p)}/\sqrt{n}.
\end{equation}
For large $p$,
$\hat{\Omega}$ is a noisy estimate for $\Omega$, so we regularize it by
%
%
\begin{equation} \label{Defineomega*}
\Omega^*(i,j) = \hat{\Omega}(i,j) 1_{\{ | \hat{\Omega}(i,j)| \geq
\log^{-1}(p) \}}.
\end{equation}
The threshold $\log^{-1}(p)$ is chosen mainly for simplicity and can
be replaced by $\log^{-a}(p)$, where $a > 0$ is a constant.
The following lemma is a direct result of (\ref{Bickel}); we omit the proof.
%
%
\begin{lemma} \label{lemmasparsifying}
Fix $A > 0$, $\gamma\in(0,1)$ and $\omega_0 \in(0, 1/2)$. As $p
\goto\infty$, for any $\Omega\in{\mathcal M}_p^*(\omega_0, \gamma
, A)$, with probability of $1 - o(1/p^2)$, each row of $\Omega^*$ has
no more than $2 \log(p)$ nonzero coordinates, and $\|\Omega^* -
\Omega\|_{\infty} \leq C (\log(p))^{-(1-\gamma)}$.
\end{lemma}

Taking $\Omega_0 = \Omega^*$, we form a graph $(V_0, \Omega^*)$.
The following lemma is proved in \cite{UPSSupp}, which says that,
except for a negligible probability, ${\mathcal U}_p(t_p^*)$ has
the SAS property.
%
%
\begin{lemma}[(SAS)] \label{lemmaSAS}
Consider model (\ref{modelrecap})--(\ref{Assumptiona})
where (\ref{Assumptionb})--(\ref{Assumptione}) hold.
Set $t^*_p$ as (\ref{Definet}).
As $p \goto\infty$, there is a constant $K$ such that
with probability $1 - L_p p^{-(\vartheta+ r)^2/(4r)}$, ${\mathcal
U}_p(t_p^*)$ has the SAS property with respect to $(V_0, \Omega^*, K)$.
\end{lemma}
%

\subsection{Reduction to many small-size regression problems}
Together, the sure screening property and the SAS property make sure that
the original regression problem reduces to many separate small-size regression
problems. In detail,
the SAS property implies that ${\mathcal U}_p(t_p^*)$ splits into many
connected subgraphs, each is small in size, and different ones are
disconnected. Given two disjoint connected subgraphs $\call$ and
${\mathcal J}_0$ where $\call\lhd\mathcal{U}_p(t)$ and ${\mathcal
J}_0 \lhd\mathcal{U}_p(t)$,
%
%
\begin{equation} \label{Disconnected}
\Omega^*(i, j) = 0\qquad \forall i \in\call, j \in
{\mathcal J}_0.
\end{equation}
Recall that the regression model (\ref{model1}) is closely related to
the model $X' Y = X'X \beta+ X'z$. Fixing a connected subgraph $\call
\lhd{\mathcal U}_p(t_p^*)$, we restrict our attention to $\call$ by
considering
$(X' Y)^{\call} = (X'X \beta)^{\call} + (X' z)^{\call}$.
See\vspace*{-1pt} Definition \ref{def1.1} for notation.
Since\vspace*{1pt} $X_i \stackrel{\mathrm{i.i.d.}}{\sim} N(0, \frac{1}{n}\Omega)$ and
$\call$ has a small size, we expect to see $(X'X \beta)^{\call}
\approx(\Omega\beta)^{\call}$ and
$(X'z)^{\call} \approx N(0, \Omega^{\call, \call})$. Therefore,
$(X' Y)^{\call} \approx N((\Omega\beta)^{\call}, \Omega^{\call,
\call})$.
A~key observation is
%
%
\begin{equation} \label{restrictionequivalence}
(\Omega\beta)^{\call} \approx\Omega^{\call, \call} \beta^{\call}.
\end{equation}
In fact, letting $\call^c = \{j\dvtx1 \leq j \leq p, j \notin\call\}$,
it is seen that
%
%
\begin{equation} \label{point}\qquad
(\Omega\beta)^{\call} - \Omega^{\call, \call} \beta^{\call} =
(\Omega^*)^{\call, \call^c} \beta^{\call^c} + (\Omega- \Omega
^*)^{\call, \call^c} \beta^{\call^c} = \mbox{I} + \mbox{II}.
\end{equation}
First, by Lemma \ref{lemmasparsifying}, $| \mbox{II} | \leq C\|\Omega-
\Omega^* \|_{\infty} \|\beta\|_{\infty} = o (\sqrt{\log(p)} )$
coordinate-wise, hence II is negligible. Second, by the sure screening\vadjust{\goodbreak}
property, signals that are falsely screened out in the $U$-step are
fewer than $L_p p^{1-(\vartheta+ r)^2/(4r)}$, and therefore have a
negligible effect. To bring out the intuition, we assume~${\mathcal
U}_p(t_p^*)$ contains all signals for a moment (see \cite{UPSSupp},
Lemma A.4, for formal treatment). This, with~(\ref{Disconnected}),
implies that $\mbox{I} = 0$, and (\ref{restrictionequivalence}) follows.

As a result, the original regression problem reduces to many small-size
regression problems of the form
%
%
\begin{equation} \label{Pmodel}
(X' Y)^{\call} \approx N( \Omega^{\call, \call} \beta^{\call},
\Omega^{\call, \call})
\end{equation}
that can be fitted separately. Note that $\Omega^{\call, \call}$ can
be accurately estimated by $(X'X)^{\call, \call}$, due to the small
size of $\call$. We are now ready for the $P$-step.

\subsection{$P$-step}
The goal of the $P$-step is that, for each fixed connected subgraph
$\call\lhd\mathcal{U}_p(t_p^*)$, we fit model
(\ref{Pmodel}) with an error rate $\leq L_p p^{-(\vartheta+
r)^2/(4r)}$. This turns out to be rather
delicate, and many methods (including the lasso and the subset
selection) do not
achieve the desired rate of convergence.

For this reason, we proposed a penalized-MLE approach. The idea can be
explained as follows.
Given that $\call\lhd{\mathcal U}_p(t_p^*)$ as a priori, the chance
that $\call$ contains $k$ signals
is \mbox{$\sim$}$\eps_p^k$. This motivates us to fit model (\ref{Pmodel}) by
maximizing the likelihood function
$\eps_p^k \cdot\operatorname{exp} [- \frac{1}{2} [( X'Y)^{\call}
- A \mu]' A^{-1} [(X'Y)^{\call} - A \mu] ]$,
subject to $\|\mu\|_0 = k$.
Recalling $A=(X'X)^{\call, \call} \approx\Omega^{\call, \call}$,
this is proportional to the density of $(X'Y)^{\call}$ in (\ref
{Pmodel}), hence the name of penalized MLE. Recalling $\eps_p =
p^{-\vartheta}$ and $\lambda_p^{\mathrm{ups}} = \sqrt{2 \vartheta\log p}$,
it is equivalent to minimizing
%
%
\begin{equation} \label{MainPstep}
[(X'Y)^{\call} - A \mu]' A^{-1} [(X'Y)^{\call
} - A \mu] + (\lambda_p^{\mathrm{ups}} )^2 \cdot\| \mu\|_0.
\end{equation}

Unfortunately, (\ref{MainPstep}) does not achieve the desired rate of
convergence as expected. The reason is that we have not taken full
advantage of the information provided: given that all coordinates in
$\call$ survive the screening, each signal in $\call$ should be
relatively strong. Motivated by this, for some tuning parameter
$u^{\mathrm{ups}} > 0$, we force all nonzero coordinates of $\mu$ to
equal~$u^{\mathrm{ups}}$. This is the UPS procedure we introduced in Section \ref
{secIntro}. In Theorem \ref{thmUB} below, we show that this procedure
obtains the desired rate of convergence provided that $u^{\mathrm{ups}}$ is
properly set.

One may think that forcing all nonzero coordinates of $\mu$ to be
equal is too restrictive, since the nonzero coordinates of $\beta
^{\call}$ are unequal. Nevertheless, the UPS achieves the desired
error rate.
The reason is that, knowing
the exact values of the nonzero coordinates is not crucial, as the main
goal is to separate nonzero coordinates of $\beta^{\call}$ from the
zero ones.

Similarly, since knowing the signal distribution $\pi_p$ may be very
helpful, one may choose to
estimate $\pi_p$ using the data first and then combine the estimated
distribution with the $P$-step.
However, this has two drawbacks. First, model (\ref{Pmodel}) is very
small in size,
and can be easily over fit if
we introduce too many degrees of freedom.
Second, estimating $\pi_p$ usually involves deconvolution, which
generally has
relatively slow rate of convergence (e.g.,~\cite{Allofstat}); a noisy
estimate of $\pi_p$ may
hurt rather than help in fitting model~(\ref{Pmodel}).

\subsection{Upper bound}
We are now ready for the upper bound. To recap,
the proposed procedure is as follows:
\begin{itemize}
\item With fixed tuning parameters $(t, \lambda^{\mathrm{ups}}, u^{\mathrm{ups}})$,
obtain ${\mathcal U}_p(t) = \{j\dvtx1 \leq j \leq p$, $(x_j, Y) \geq t\}$.
\item Obtain $\Omega^*$ as in (\ref{Defineomega*}), and form a graph
$(V_0, \Omega_0)$ with $V_0 = \{1, 2, \ldots, p\}$ and $\Omega_0 =
\Omega^*$.
\item Split ${\mathcal U}_p(t)$ into connected subgraphs where
different ones are disconnected. For each connected subgraph $\call=
\{ i_1,i_2, \ldots, i_K\}$, obtain the minimizer of (\ref
{MainPstep}), where each coordinate of $\mu$ is either $0$ or
$u^{\mathrm{ups}}$. Denote the estimate by $\hat{\mu}(\call) = \hat{\mu
}(\call; Y, X, t, \lambda^{\mathrm{ups}}, u^{\mathrm{ups}}, p)$.
\item For any $1 \leq j \leq p$, if $j \notin{\mathcal U}_p(t)$, set
$\hb_j = 0$. Otherwise, there is a unique $\call= \{ i_1, i_2, \ldots
, i_K \} \lhd{\mathcal U}_p(t)$, where $i_1 < i_2 <\cdots< i_K$,
such that $j$ is the $k$th coordinate of $\call$. Set $\hb_j = (\hat
{\mu}(\call))_k$.
\end{itemize}
Denote the resulting estimator by $\hb(Y, X; t, \lambda^{\mathrm{ups}},
u^{\mathrm{ups}})$. We have the following theorem.
%
%
\begin{theorem} \label{thmUB}
Consider model (\ref{modelrecap})--(\ref{Assumptiona}) where (\ref
{Assumptionb})--(\ref{Assumptione}) hold, and fix $0 < q \leq
(\vartheta+ r)^2/(4r)$. For sufficiently large $p$, if $\Omega^{(p)}
\in{\mathcal M}_p^+(\omega_0, \gamma, A)$, and we
set the tuning parameters of the UPS at
\[
t = t_p^* = \sqrt{2 q \log(p)},\qquad \lambda^{\mathrm{ups}} = \lambda
_p^{\mathrm{ups}} = \sqrt{2 \vartheta\log p},\qquad u^{\mathrm{ups}} = u_p^{\mathrm{ups}} =
\tau_p,
\]
then as $p \goto\infty$, $\hamm_p(\hb^{\mathrm{ups}}(Y, X; t_p^*,
\lambda_p^{\mathrm{ups}}, u_p^{\mathrm{ups}}), \vartheta, r,
\Omega^{(p)}) \leq L_p \cdot s_p \cdot\break p^{- {(r -
\vartheta)^2}/({4r})}$. The claim remains valid if $r/\vartheta\leq3 + 2
\sqrt{2}$ and $\Omega^{(p)} \in{\mathcal M}_p^*(\omega_0,\break \gamma, A)$
for sufficiently large $p$.
\end{theorem}

Except for the $L_p$ term, the upper bound matches the lower bound in
Theorem~\ref{thmLB}. Therefore, both bounds are tight and the UPS is
rate optimal.

\subsection{Tuning parameters of the UPS} \label{secTuning}
The UPS uses three tuning parameters $(t_p^*, \lambda_p^{\mathrm{ups}},
u_p^{\mathrm{ups}})$. In this section, we show that
under certain conditions, the parameters $(\lambda_p^{\mathrm{ups}},
u_p^{\mathrm{ups}})$ can be estimated from the data.

In detail, recall that $\tilde{Y} = X' Y$. For $t > 0$, introduce
$\bar{F}_p(t) = \frac{1}{p} \sum_{j = 1}^p 1 \{ \ty_j > t\} $ and $
\mu_p(t) = \frac{1}{p} \sum_{j = 1}^p \ty_j \cdot1\{ \ty_j > t\}$.
Denote the largest off-diagonal coordinate of $\Omega$ by
$\delta_0 = \delta_0(\Omega) = {\max_{\{ 1 \leq i, j \leq p, i \neq
j \}}} |\Omega(i,j)|$.
Recalling that the support of $\pi_p$ is contained in $[\tau_p,
(1+\eta)\tau_p]$, we suppose
%
%
\begin{equation} \label{delta0condition}
2 \delta_0 ( 1 + \eta) - 1 \leq\vartheta/r\qquad \mbox{so that }
\delta_0^2 (1+\eta)^2r < \frac{(\vartheta+ r)^2}{4r}.
\end{equation}
Let $\mu_p^*(\pi_p)$ be the mean of $\pi_p$.
The following is proved in \cite{UPSSupp}.
%
%
\begin{lemma} \label{lemmatuningstoch}
Fix $q$ such that $\max\{\delta_0^2 (1+\eta)^2r, \vartheta\} < q
\leq(\vartheta+ r)^2/(4r)$, and let $t_p^* = \sqrt{2 q \log p}$.
Suppose\vspace*{1pt} the conditions in Theorem \ref{thmUB} hold. As \mbox{$p \goto\infty$},
with probability of $1 - o(1/p)$,
%
%
\begin{equation}\label{largedevcondition}
|[\cf_p(t_p^*) / \eps_p] - 1 | = o(1)\quad
\mbox{and}\quad
|[ \mu_p(t_p^*) / (\varepsilon_p \mu_p^*(\pi_p))] - 1| = o(1).\hspace*{-30pt}
\end{equation}
\end{lemma}

Motivated by Lemma \ref{largedevcondition}, we propose to estimate
$(\lambda^{\mathrm{ups}}, u^{\mathrm{ups}})$ by
%
%
\begin{eqnarray} \label{UPStuningAdd}
\hat{\lambda}_p^{\mathrm{ups}} &=& \hat{\lambda}_p^{\mathrm{ups}}(q) = \sqrt{-2 \log
(\cf_p(t_p^*))},\nonumber\\[-8pt]\\[-8pt]
\hat{u}_p^{\mathrm{ups}} &=& \hat{u}_p^{\mathrm{ups}}(q) = \mu
_p(t_p^*)/\cf_p(t_p^*).\nonumber
\end{eqnarray}

%
%
\begin{theorem} \label{thmadaptive}
Fix $q$ such that $\max\{\delta_0^2 (1+\eta)^2r, \vartheta\} < q
\leq(\vartheta+ r)^2/(4r)$, and let $t_p^* = \sqrt{2 q \log p}$.
Suppose the conditions of Theorem \ref{thmUB} hold. As $p \goto\infty
$, if additionally
$\mu_p^*(\pi_p) \leq(1 + o(1)) \tau_p$,
then
$\hamm_p(\hb^{\mathrm{ups}}) \leq
L_p \cdot s_p \cdot p^{- (r - \vartheta)^2/(4r)}$.
\end{theorem}

As a result, $t_p^*$ is the only tuning parameter needed by the UPS.
By Theorem \ref{thmadaptive}, the
performance of the UPS is relatively insensitive to the choice of
$t_p^*$, as long as it falls in a certain range. Numerical studies in
Section \ref{secSimul} confirm this for finite $p$. The numerical
study also suggests that the lasso is comparably more sensitive to its
tuning parameter $\lambda^{\mathrm{lasso}}$.

\subsection{Discussions}
While the conditions in Theorems \ref{thmUB} are \ref{thmadaptive} are
relatively strong, the key idea of the
paper applies to much broader settings.
The success of UPS attributes to the interaction of the signal sparsity
and graph sparsity, which can be found in many
applications [e.g., compressive sensing, genome-wide association study (GWAS)].

In the forthcoming papers \cite{FJZ,JinZhang1,JinZhang2}, we revisit
the key idea of this paper, and extend our results to more general
settings. However, the current paper is different from
\cite{FJZ,JinZhang1,JinZhang2} in important ways.
First, the focus of \cite{FJZ} is on ill-posed regression models
and change-point problems, and the focus of~\cite{JinZhang2} is on
Ising model and network data. Second, the current paper uses the
so-called ``phase diagram'' as a new criterion for optimality (e.g.,
\cite{DonohoTanner}), and Jin and Zhang \cite{JinZhang1} use the more
traditional ``asymptotic minimaxity'' as the criterion for optimality.
Due to the complexity of the problem, one type of optimality usually
does not imply the other. The current paper and \cite{JinZhang1} have
very different targets, objectives and underlying mathematical
techniques, and the results in either one cannot be deduced from the
other.

The current paper is new in at least two aspects.
First, given that marginal regression is a widely used method but is
not well justified, this paper shows that marginal regression can
actually work, provided that
an additional cleaning stage is performed. Second, it shows that
$L^0$-penalization method---the target of many relaxation methods---is
nonoptimal,
even in very simple settings, and even when the tuning parameter is
ideally set.

\section{A refinement for moderately large $p$} \label{secrefinement}

We introduce a refinement for the UPS when $p$ is moderately large.
We begin by investigating the relationship between the regression model
and Stein's normal means model.

Recall that model (\ref{model1}) is closely related to the following model:
%
%
\begin{equation} \label{XY}
X' Y =X' X \beta+ X' z,\qquad z \sim N(0, I_n),
\end{equation}
which is approximately equivalent to Stein's normal means model as follows:
%
%
\begin{equation} \label{XY2}
X' Y \approx\Omega\beta+ N(0, \Omega) \quad\Longleftrightarrow\quad
\Omega^{-1} X'Y \approx N(\beta, \Omega^{-1}).
\end{equation}
In the literature, Stein's normal means model has been extensively
studied, but the focus has been on the case where $\Omega$ is
diagonal (e.g., \cite{Allofstat}). When~$\Omega$ is not diagonal,
Stein's normal means model is intrinsically a regression problem. To
see how close models (\ref{XY}) and (\ref{XY2}) are, write
%
%
\begin{equation} \label{XY3}
X' Y\,{=}\,\biggl[ \Omega\beta\,{+}\,\frac{\sqrt{n}}{\|z\|} X'z
\biggr]\,{+}\,\biggl[ (X' X\,{-}\,\Omega) \beta\,{+}\,\biggl(\frac{\|z\|}{\sqrt{n}}\,{-}\,1
\biggr) \frac{\sqrt{n}}{\|z\|} X'z \biggr]\,{=}\,\mbox{I}\,{+}\,\mbox{II}.\hspace*{-35pt}
\end{equation}
First, note that $\mbox{I} \sim N(\Omega\beta, \Omega)$. For $\mbox{II}$, we have
the following lemma. 
%
%
\begin{lemma} \label{lemmaStochmean}
Consider model (\ref{modelrecap})--(\ref{Assumptiona}) where (\ref
{Assumptiona})--(\ref{Assumptionc}) hold. As\break \mbox{$p \goto\infty$}, there
is a constant $C > 0$ such that
except for a probability of $o(1/p)$,
\begin{eqnarray*}
\biggl|\frac{\|z\|}{\sqrt{n}} -1 \biggr| &\leq& C \bigl(\sqrt{\log p}\bigr)
p^{-\theta/2},\\
\|(X'X - \Omega) \beta\|_{\infty} &\leq& C \|
\Omega\| \bigl(\sqrt{2\log p}\bigr) p^{-({\theta- (1 - \vartheta)})/{2}}.
\end{eqnarray*}
\end{lemma}

It follows that $| \mbox{II} | \leq C \sqrt{2\log(p}) \cdot p^{- [\theta-
(1 - \vartheta)]/2}$ coordinate-wise. Therefore, \textit{asymptotically},
models (\ref{XY}) and (\ref{XY2}) have negligible difference.
However, when $p$ is moderately large, the difference between models
(\ref{XY}) and~(\ref{XY2}) may be nonnegligible. In Table \ref
{tabDesignError},
we tabulate the values of $\sqrt{2\log(p}) \cdot p^{- [\theta- (1 -
\vartheta)]/2}$, which are relatively large for moderately large $p$.

%
%
\begin{table}[b]
\caption{The values of $\sqrt{2 \log(p)} p^{-[\theta- (1 -
\vartheta)]/2} $ for different $p$ and $(\theta, \vartheta)$}
\label{tabDesignError}
\begin{tabular*}{\tablewidth}{@{\extracolsep{\fill}}lcccccc@{}}
\hline
$\bolds{p}$ & $\bolds{400}$ & $\bolds{5 \times400}$ &
$\bolds{5^2 \times400}$ & $\bolds{5^3 \times400}$ & $\bolds{5^4
\times400}$ & $\bolds{5^5 \times400}$\\
\hline
$(\theta, \vartheta)=(0.91, 0.65)$ & $0.65$ & $0.46$ & $0.33$
&$0.22$& $0.15$& $0.10$ \\
$(\theta, \vartheta)=(0.91, 0.5)$
& $1.01$ & $0.82$ &$0.65$ &$0.51$& $0.39$ &$0.30$ \\
\hline
\end{tabular*}
\end{table}

This says that, for moderately large $p$, the random design model is
much noisier than Stein's normal means model. As a result, in the
$U$-step, we tend to falsely keep more noise terms in the former than
in the latter; some of these\vadjust{\goodbreak} noise terms are large in magnitude, and it
is hard to clean all of them in the $P$-step. To see how the problem
can be fixed, we write
%
%
\begin{equation} \label{troubleterm}
X'X \beta= (X'X - \Omega^*) \beta+ \Omega^* \beta.
\end{equation}
On one hand, the term $(X'X - \Omega^*) \beta$ causes the random
design model to be much noisier than Stein's normal means model. On the
other hand, this term can be easily removed from the model if we have a
reasonably good estimate of $\beta$. This motivates a refinement as follows.

For any $p\times1$ vector $y$, let $S^2(y) = \frac{1}{p-1} \sum_{j =
1}^p (y_j - \bar{y})^2$ where $\bar{y} = \frac{1}{p}\sum_{j =1}^p
y_j$. We propose the following procedure:
(1) Run the UPS and obtain an estimate of $\beta$, say, $\hat{\beta
}$. Let $W^{(0)}=X'Y$ and $\hat{\beta}^{(0)}=\hat{\beta}$.
(2) For $j=1,2,3$, respectively, let $W^{(j)} = X'Y - (X' X - \Omega
^*) \hat{\beta}^{(j-1)}$. If $S(W^{(j)})/ S(W^{(j-1)}) \leq1.05$,
run the UPS with $X'Y$ replaced by $W^{(j)}$ and other parts unchanged,
and let $\hat{\beta}^{(j)}$ be the new estimate. Stop otherwise.

Numerical studies in Section \ref{secSimul} suggest that the
refinement is beneficial for moderately large $p$. When $p$ is
sufficiently large [e.g.,
$\sqrt{2 \log(p)} \cdot p^{-[\theta- (1 - \vartheta)]/2} \leq
0.4$], the original UPS is usually good enough. In this case,
refinements are not necessary, but may still offer improvements.

\section{Understanding the lasso and the subset selection} \label{seclasso}

In this section,
we show that there is a region in the phase space where the lasso is
rate nonoptimal (similarly for subset selection).
We use Stein's normal means
model instead of the random design model
(as the goal is to understand the nonoptimality
of these methods, focusing on a simpler model enjoys mathematical
convenience, yet is also sufficient; see Section \ref{subseclasso}).

To recap, the model we consider in this section is
$\tilde{Y} \sim N(\Omega\beta, \Omega)$,
where $\tilde{Y}$ is the counterpart of $X'Y$ in the random design
model. Fix $a \in(-1/2, 1/2)$. As in
Section \ref{subseclasso}, we let $\Omega$ be the tridiagonal matrix
as in (\ref{tridiagonal2}), and
$\pi_p$ be the point mass at $\tau_p = \sqrt{2 r \log p}$. In other words,
%
%
\begin{equation} \label{stein1b}
\beta_j \stackrel{\mathrm{i.i.d.}}{\sim} (1 - \eps_p) \nu_0 + \eps_p \nu
_{\tau_p},\qquad \eps_p = p^{-\vartheta}, \qquad \tau_p = \sqrt
{2 r \log p}.
\end{equation}
Throughout this section, we assume $r > \vartheta$ so that successful
variable selection is possible. Somewhat surprisingly, even in this
simple case and even when $(\eps_p, \tau_p)$ are known, there is a
region in the phase space where neither the lasso nor the subset
selection is optimal. To shed light, we first take a~heuristic approach
below. Formal statements are given later.

\subsection{Understanding the lasso}
The vector $\ty$ consists of three main components: true signals, fake
signals and pure noise (see Definition \ref{def1.3}).
According to~(\ref{stein1b}), true signals may appear
as singletons, pairs, triplets, etc., but singletons are the most
common and therefore have the major effect.
For each signal singleton, since $\Omega$ is tridiagonal, we have two
fake signals, one to the left and one to the right. Given a site $j$,
$1 \leq j \leq p$,
the lasso may make three types of errors:
\begin{itemize}
\item\textit{Type} I. $\ty_j$ is a pure noise, but the lasso mistakes it
as a signal.\vspace*{1pt}
\item\textit{Type} II. $\ty_j$ is a signal singleton, but the lasso
mistakes it as a noise.\vspace*{1pt}
\item\textit{Type} III. $\ty_j$ is a fake signal next to a signal
singleton, but the lasso mistakes it as a signal.
\end{itemize}
There are other types of errors, but these are the major ones.

To minimize the sum of these errors, the lasso needs to choose the
tuning parameter $\lambda^{\mathrm{lasso}}$ carefully. To shed light, we
first consider the uncorrelated case where $\Omega$ is the identity
matrix. In this case, we do not have fake signals and it is
understood that the lasso is equivalent to the soft-thresholding
procedure \cite{Allofstat}, where the expected sum of types I and
II errors is
%
%
\begin{equation} \label{lassohamm0}
p [ (1 - \eps_p) \bphi(\lambda^{\mathrm{lasso}}) + \eps_p \Phi(
\lambda^{\mathrm{lasso}}-\tau_p) ].
\end{equation}
Here, $\bphi= 1 - \Phi$ is the survival function of $N(0,1)$.
In (\ref{lassohamm0}), fixing $0 < q < 1$ and taking
$\lambda^{\mathrm{lasso}} = \lambda_p^{\mathrm{lasso}} = \sqrt{2 q \log(p)}$,
the expected sum of errors is
\[
\sim\cases{
L_p \bigl[p^{1-q} + p^{1-(\vartheta+ (\sqrt{q} - \sqrt{r})^2 )}
\bigr], &\quad if $0<q<r$,
\vspace*{2pt}\cr
p^{1-q} + p^{1-\vartheta},&\quad if $q> r$.}
\]
The right-hand\vspace*{-1pt} side is minimized at $q = (\vartheta+ r)^2/(4r)$ at which
$\lambda_p^{\mathrm{lasso}} = \frac{\vartheta+ r}{2r} \tau_p$, and the sum
of errors is $L_p p^{1-(\vartheta+ r)^2/(4r)}$, which is the optimal
rate of convergence.
For a smaller $q$, the lasso keeps too many noise terms. For a larger
$q$, the lasso kills too many signals.

Return to the correlated case. The vector $\tilde{Y}$ is at least as
noisy as that in the uncorrelated case. As a result, to
control the type I errors, we should choose $\lambda_p^{\mathrm{lasso}} $ to be
at least
$\frac{\vartheta+ r}{2r} \tau_p$.
This is confirmed in Lemma \ref{lemmalasso} below.

In light of this, we fix $q \geq(\vartheta+ r)^2/(4r)$ and let
$\lambda_p^{\mathrm{lasso}} = \sqrt{2 q \log(p)}$ from now on. We observe that
except for a negligible probability, the support of $\hb^{\mathrm{lasso}}$,
denoted by $\hat{S}_p^{\mathrm{lasso}}$, splits into many small clusters (i.e.,
block of adjacent indices). There is an integer $K$ not depending on
$p$ that has the following effects: (a)
If $\tilde{Y}_j$ is a pure noise, and there is no signal within\vspace*{1pt}
a~distance of $K$ from it, then
either $\hb_j^{\mathrm{lasso}} = 0$, or $\hb_j^{\mathrm{lasso}} \neq0$ but $\hb_{j
\pm1}^{\mathrm{lasso}} = 0$, and (b)
If $\tilde{Y}_j$ is a signal singleton, and there is no other signal
within a~distance of $K$ from it, then
either $\hb_j^{\mathrm{lasso}} = 0$, or $\hb_j^{\mathrm{lasso}} \neq0$ but $\hb_{j
\pm2} = 0$ and at least one of $\{ \hb_{j + 1}^{\mathrm{lasso}}, \hb
_{j-1}^{\mathrm{lasso}}\} $ is $0$.
These\vspace*{1pt} heuristics are justified in \cite{Thesis}
(we use such heuristics to provide insight, but not for proving results below).

\begin{figure}

\includegraphics{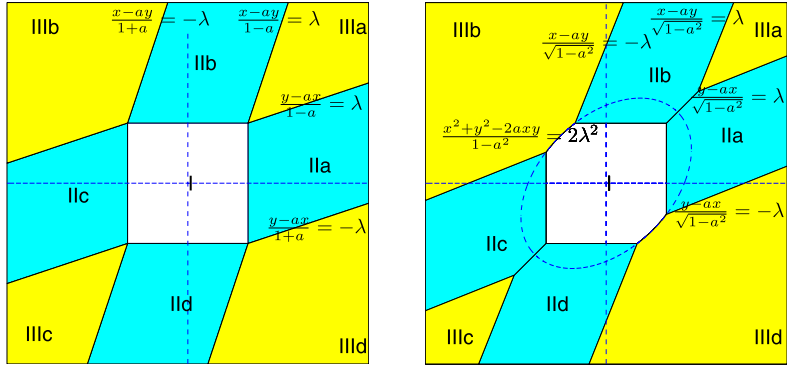}

\caption{Partition of regions as in Lemma \protect\ref{lemma2dlasso} (left)
and in Lemma \protect\ref{lemma2dss} (right).}
\label{fig2dlasso}
\end{figure}

At the same time, let $\call= \{j-k+1, \ldots, j \} \subset\hat
{S}_p^{\mathrm{lasso}}$ be a cluster, so that
$\hb_{j-k}^{\mathrm{lasso}} = \hb_{j + 1} ^{\mathrm{lasso}} = 0$.
Since\vspace*{2pt} $\Omega$ is tridiagonal,
$(\hb^{\mathrm{lasso}})^{\call}$, the restriction of $\hb
^{\mathrm{lasso}}$ to $\call$, is the solution of the following
small-size minimization problem:
%
%
\begin{equation} \label{2dlasso1}
\tfrac{1}{2} \mu' (\Omega^{\call, \call}) \mu- \mu' \tilde
{Y}^{\call} + \lambda^{\mathrm{lasso}} \| \mu\|_1\qquad \mbox{where $\mu$
is a $k\times1$ vector}.\hspace*{-30pt}
\end{equation}
See Definition \ref{def1.1}. Two special cases are noteworthy. First, $\call=
\{ j \}$, and the solution of (\ref{2dlasso1}) is given by
$\hb_j^{\mathrm{lasso}} = \operatorname{sgn}( \tilde{Y}_j) (|\tilde{Y}_j| -
\lambda^{\mathrm{lasso}})^+$,
which is the soft-thresholding \cite{Allofstat}. Second, $\call= \{j
-1, j\}$. We call the solution of (\ref{2dlasso1}) in this case the
\textit{bivariate lasso}. We have the following lemma, where all regions
I-IIId are illustrated in Figure \ref{fig2dlasso} ($x$-axis is
$\ty_{j-1}$, $y$-axis is $\ty_j$).
%
%
\begin{lemma} \label{lemma2dlasso} Denote $\lambda= \lambda
^{\mathrm{lasso}}$. The solution of the bivariate lasso
$(\hb_{j-1}^{\mathrm{lasso}}$,
$\hb_{j}^{\mathrm{lasso}})$ is given by $(\hb_{j-1}^{\mathrm{lasso}}, \hb_{j}^{\mathrm{lasso}})
= (\operatorname{sgn}( \tilde{Y}_{j-1}) (|\tilde{Y}_{j-1}| - \lambda)^+,
\operatorname{sgn}( \tilde{Y}_j) (|\tilde{Y}_j| - \lambda)^+)$ if
$(\ty_{j-1}, \ty_j)$ is in regions \textup{I}, \textup{IIa-IId} and
$(\hb_{j-1}^{\mathrm{lasso}}, \hb_{j}^{\mathrm{lasso}}) = \frac{1}{1-a^2}(Z_{j-1} - a
Z_j, Z_j - a Z_{j -1})$
if $(\ty_{j-1}, \ty_j)$ is in regions \textup{IIIa}-\textup{IIId}. Here, $Z_{j-1} =
\ty_{j-1} - \lambda$ if $(\ty_{j-1}, \ty_j)$ is in regions \textup{IIIa},
\textup{IIId} and $Z_{j-1} = \ty_{j-1} + \lambda$ otherwise; $Z_{j} = \ty
_{j} - \lambda$ if $(\ty_{j-1}, \ty_j)$ is in regions \textup{IIIa, IIIb} and
$Z_j = \ty_j + \lambda$ otherwise.
\end{lemma}

In the white region of Figure \ref{fig2dlasso}, both
$\hb_{j-1}^{\mathrm{lasso}}$ and $\hb_j^{\mathrm{lasso}}$ are $0$. In the blue regions,
exactly one of them is $0$. In the yellow regions, both are nonzero.
Lemma \ref{lemma2dlasso} is proved in \cite{UPSSupp}.

%
%

As a result, the following hold, except for a negligible probability:
\begin{itemize}
\item\textit{Type} I. There are $O(p)$ indices $j$ where $\tilde{Y}_j$
is a pure noise, and no signal appears within a distance of $K$ from
it. For each of such $j$, the lasso acts on $\ty_j$ as (univariate)
soft-thresholding, and $\hb_j^{\mathrm{lasso}} \neq0$ if and only if $|\ty_j|
\geq\lambda_p^{\mathrm{lasso}} $.
\item\textit{Types} II--III. There are $O(p \eps_p)$ indices where $\ty
_j$ is a signal singleton, and no other signal appears within a
distance of $K$ from it. The lasso either acts on $\ty_j$ as
soft-thresholding, or acts on both $\ty_j$ and one of its neighbors as
the bivariate lasso. As a result, $\hb_j^{\mathrm{lasso}} = 0$ if and only if
$|\ty_j| \leq\lambda_p^{\mathrm{lasso}} $ (type~II), and both $\hb
_j^{\mathrm{lasso}}$ and $\hb_{j-1}^{\mathrm{lasso}}$ are nonzero if and only if $(\ty
_{j-1}, \ty_j)'$ falls in regions IIIa-IIId, with IIIa and
IIIb being the most likely (type~III).\vadjust{\goodbreak}
\end{itemize}
Noting that $\ty_j \sim N(0,1)$ if it is a pure noise and $\ty_j \sim
N(\tau_p, 1)$ if it is a~signal singleton, the sum of types I and
II errors is
$L_p p [ P( N(0,1) \geq\lambda_p^{\mathrm{lasso}} ) + \eps_p P (N(\tau
_p, 1) < \lambda_p^{\mathrm{lasso}} ) ] = L_p p [\bphi(\lambda
_p^{\mathrm{lasso}} ) + \eps_p \Phi(\lambda_p^{\mathrm{lasso}} - \tau_p) ]$.
Also, when $\ty_j$ is a~signal singleton, $(\ty_{j-1}, \ty_j)'$ is
distributed as a~bivariate normal with\break means~$a \tau_p$ and $\tau_p$,
variances $1$, and correlation $a$. Denote such a bivariate normal
distribution by $W$ for short. The type III error is
\mbox{$L_p p \cdot P(\beta_{j-1} = 0$}, \mbox{$\beta_j = \tau_p$}, $(\ty
_{j-1}, \ty_j)' \in\mbox{regions IIIa or IIIb} ) \sim L_p
p \eps_p \cdot P(W \in\mbox{regions IIIa}\break \mbox{or IIIb} )$.
Therefore, the sum of three types of errors is
%
%
\begin{equation} \label{lassohamm}
L_p p \cdot[ \bphi(\lambda_p^{\mathrm{lasso}} ) + \eps_p \Phi
(\lambda_p^{\mathrm{lasso}} - \tau_p) + \eps_p P(W \in\mbox{regions
IIIa or IIIb} ) ],\hspace*{-30pt}
\end{equation}
which can be conveniently evaluated.
Note that the sum of types I and II errors in the correlated
case is the same as that in the uncorrelated case, which
is minimized at $\lambda_p^{\mathrm{lasso}} = (\vartheta+ r)/(2r) \tau_p$.
Therefore, whether the lasso is optimal or not depends on whether the
type III error is smaller than the optimal rate of convergence or not.
Unfortunately, in certain regions of the phase space, the type III
error can be significantly larger than the optimal rate.
In other words, provided that the tuning parameters are properly set,
the lasso is able to separate the signal singletons from the pure noise.
However, it may not be efficient in filtering out the fake signals,
which is the culprit for its nonoptimality.

For short, write $\hamm_p(\hat{\beta}^{\mathrm{lasso}}(\lambda_p^{\mathrm{lasso}})) =
\hamm(\hat{\beta}^{\mathrm{lasso}}(\lambda_p^{\mathrm{lasso}}); \eps_p, \tau_p, a)$.
The following is proved in \cite{UPSSupp}, confirming the above
heuristics.
%
%
\begin{lemma} \label{lemmalasso}
Fix $\vartheta\in(0,1)$, $r > \vartheta$, $q > 0$ and $a \in
(-1/2,1/2)$. Set the lasso tuning parameter as $\lambda^{\mathrm{lasso}}_p =
\sqrt{2 q \log p}$. As $p \goto\infty$,
\begin{eqnarray*}
&&\frac{\hamm(\hat{\beta}^{\mathrm{lasso}}(\lambda_p^{\mathrm{lasso}}))}{s_p} \\
&&\qquad\geq
\cases{
L_p p^{-\min\{(({1-|a|})/({1+|a|}))q, q-\vartheta\}}, &\quad if $0 < q <
\dfrac{(\vartheta+ r)^2}{4r}$,
\vspace*{2pt}\cr
L_p p^{-\min\{(({1-|a|})/({1+|a|}))q, (\sqrt{r}-\sqrt{q})^2 \}}, &\quad
if $\dfrac{(\vartheta+ r)^2}{4r} < q < r$,
\vspace*{2pt}\cr
\bigl(1+o(1)\bigr), &\quad if $q > r$.}
\end{eqnarray*}
\end{lemma}

The exponent on the right-hand side
is minimized at
$q = (\vartheta+ r)^2/(4r)$ when $r < [(1 + \sqrt{1 - a^2}) / |a| ]
\vartheta$ and
$q = (1 + |a|) (1 - \sqrt{1 - a^2}) r/ (2a^2)$ when $r > [(1 + \sqrt
{1 - a^2})/|a| ] \vartheta$,
where we note that $r < [(1 + \sqrt{1 - a^2})/|a| ] \vartheta$ and $r
> [(1 + \sqrt{1 - a^2})/ |a| ] \vartheta$ correspond to the optimal
and nonoptimal regions of the lasso, respectively.
This shows that in the optimal region of the lasso, $\lambda_p^{\mathrm{lasso}}
= (\vartheta+ r)/(2r) \tau_p$ remains the optimal tuning parameter,
at which the sum of types I and II errors is minimized, and the
type III error has a negligible effect. In the nonoptimal region of
the lasso, at $\lambda_p^{\mathrm{lasso}} = (\vartheta+ r)/(2r) \tau_p$,
the type III error is larger than the sum of types I and II errors,
so the lasso needs to raise the tuning parameter slightly to
minimize the sum of all three types of errors (but the resultant
Hamming error is still larger than that of the optimal procedure).
Combining this with
Lemma \ref{lemmalasso} gives the following theorem, the proof of which
is omitted.
%
%
\begin{theorem} \label{thmlasso}
Set $\lambda^{\mathrm{lasso}}_p = \sqrt{2 q \log p}$. For all choices of $q >
0$, the error rate of the lasso satisfies $\hamm_p(\hb^{\mathrm{lasso}}(
\lambda_p^{\mathrm{lasso}})) \geq L_p \cdot s_p \cdot p^{- {(\vartheta
- r)^2}/({4r})}$ when $r / \vartheta< (1 + \sqrt{1 - a^2})/|a|$ and
\[
\hamm_p(\hb^{\mathrm{lasso}}( \lambda_p^{\mathrm{lasso}})) \geq L_p \cdot s_p \cdot
p^{\vartheta- ({(1 - |a|) (1 - \sqrt{1 - a^2}) }/({2 a^2})) r},
\]
when
$ r / \vartheta> (1 + \sqrt{1 - a^2})/|a|$.
\end{theorem}

In \cite{Thesis}, we show that when $r/\vartheta\leq3 + 2 \sqrt
{2}$, the lower bound in Theorem~\ref{thmlasso} is tight.
The proofs are relatively long, so
we leave the details to \cite{Thesis}.

\subsection{Understanding subset selection} \label{secSS}
The discussion is similar, so we keep it brief.
Fix $1 \leq j \leq p$. The major errors that subset selection makes are
the following (type III is defined differently from that in the
preceding section):
\begin{itemize}
\item
\textit{Type} I. $\ty_j$ is a pure noise, but subset selection takes it
as a signal.\vspace*{1pt}
\item
\textit{Type} II. $\ty_j$ is a signal singleton, but subset selection
takes it as a noise.
\item
\textit{Type} III. ($\ty_{j-1}, \ty_j)$ is a signal pair, but subset
selection mistakes one of them as a noise.
\end{itemize}

Suppose that $\ty_j$ is either a pure noise or a signal singleton, and
for an appropriately large $K$, no other signal appears within a
distance of $K$ from it. In this case, except for a negligible probability,
$\hb_{j \pm1}^{\mathrm{lasso}} = 0$, and the subset selection acts on site $j$
as hard thresholding \cite{Allofstat}, $\hb_j^{\mathrm{ss}} = \tilde{Y}_j
\cdot
1 \{ |\tilde{Y}_j| \geq\lambda^{\mathrm{ss}}\}$. Recall that $\tilde{Y}_j
\sim
N(0,1)$ if it is a pure noise, and $\tilde{Y}_j \sim N(\tau_p, 1)$ if
it is a signal singleton. Take $\lambda^{\mathrm{ss}} = \lambda^{\mathrm{ss}}_p =
\sqrt{2q \log p}$ as before. Similarly, the expected sum of types I
and II errors is
%
%
\begin{eqnarray} \label{ssadd1}
&&
L_p p [ \bphi(\lambda_p^{\mathrm{ss}} ) + p^{- \vartheta} \Phi(\lambda
_p^{\mathrm{ss}} - \tau_p) ] \nonumber\\[-8pt]\\[-8pt]
&&\qquad=
\cases{
L_p \bigl(p^{1-q} + p^{1- \vartheta- (\sqrt{q} - \sqrt{r})^2}\bigr), &\quad if $0 <
q < r$,
\vspace*{2pt}\cr
L_p (p^{1-q} + p^{1- \vartheta}), &\quad if $q >
r$.}\nonumber
\end{eqnarray}
On the right-hand side,
the exponent is minimized at $q = (\vartheta+ r)^2/4r$, at which the
rate is $L_p p^{1-(\vartheta+ r)^2/(4r)}$, which is the optimal rate
of convergence.

Next, consider the type III error. Suppose $(\ty_{j-1}, \ty_j)$ is a
signal pair and no other signal appears within a distance of $K$ for a
properly large $K$. Similarly, since $\Omega$ is tridiagonal,
$(\hb_{j-1}^{\mathrm{ss}}, \hb_j^{\mathrm{ss}})'$ is the\vspace*{-3pt}
minimizer of the functional $\frac{1}{2}\beta_{j-1}^2 + \frac{1}{2}
\beta_{j}^2 + a \beta_{j-1} \beta_{j} - (\tilde{Y}_{j-1} \beta_{j-1} +
\tilde{Y}_{j} \beta_j) +
\frac{(\lambda^{\mathrm{ss}}_p)^2}{2}(I\{\beta_{j-1} \neq0\} +
I\{\beta_{j}\vspace*{1pt} \neq0\})$. We call the resultant procedure
\textit{bivariate subset selection}. The following lemma
is proved in~\cite{UPSSupp}, with the regions illustrated in Figure
\ref{fig2dlasso}.
%
%
\begin{lemma} \label{lemma2dss}
The solution of the bivariate subset selection is given by
$(\hb_{j-1}^{\mathrm{ss}}, \hb_j^{\mathrm{ss}})\,{=}\,(0,0)$ if $(\tilde{Y}_{j-1},\tilde
{Y}_j)$ is in region \textup{I}, $(\hb_{j-1}^{\mathrm{ss}}, \hb_j^{\mathrm{ss}})\,{=}\,(\ty
_{j-1},0)$ if $(\tilde{Y}_{j-1},\tilde{Y}_j)$ is in regions \textup{IIa,
IIc},
$(\hb_{j-1}^{\mathrm{ss}}, \hb_j^{\mathrm{ss}}) = (0,\tilde{Y}_j)$ if $(\tilde
{Y}_{j-1},\tilde{Y}_j)$ is in regions~\textup{IIb},~i.i.d. and $(\hb
_{j-1}^{\mathrm{ss}}, \hb_j^{\mathrm{ss}}) =
(\frac{\tilde{Y}_{j-1} - a \tilde{Y}_j}{1- a ^2}, \frac{\tilde
{Y}_j - a \tilde{Y}_{j-1}}{1- a ^2})$ if $(\tilde{Y}_{j-1},\tilde
{Y}_j)$ is in regions \textup{IIIa-IIId}.\looseness=-1
\end{lemma}

When $(\tilde{Y}_{j-1}, \tilde{Y}_j)$ falls in regions I, IIa or
IIb, either $\hb_{j-1}^{\mathrm{ss}}$ or $\hb_{j}^{\mathrm{ss}}$ is $0$, and the
subset selection makes a type III error. Note there are $O(p
\eps_p^2)$ signal pairs, and that $(\tilde{Y}_{j-1}, \tilde{Y}_j)'$ is
jointly distributed as a bivariate normal with means $(1 + a) \tau_p$,
variances $1$ and correlation $a$. The type III error is then $L_p
p^{1-(2\vartheta+ \min\{[(\sqrt{r(1-a^2)}-\sqrt{q})^+]^2,
2[(\sqrt{r(1+a)}-\sqrt{q})^+]^2\} }$. Combining with (\ref{ssadd1})
and Mills's ratio gives the sum of all three types of errors.
Formally, writing for short $\hamm_p(\hb^{\mathrm{ss}}( \lambda_p^{\mathrm{ss}})) =
\hamm_p(\hb^{\mathrm{ss}}( \lambda_p^{\mathrm{ss}}); \eps_p, \tau_p, a)$, we have the
following lemma proved in \cite{UPSSupp}.
%
%
\begin{lemma} \label{lemmasubset}
Set the tuning parameter $\lambda^{\mathrm{ss}}_p = \sqrt{2 q \log p}$. The
Hamming error for the subset selection $\hamm_p(\hb^{\mathrm{ss}}( \lambda
_p^{\mathrm{ss}}))$ is at least
\[
\cases{
L_p \cdot s_p \cdot p^{-\min\{q-\vartheta, \vartheta+[(\sqrt
{r(1-a^2)}-\sqrt{q})^+]^2 \}},
&\quad if $ 0<q< \dfrac{(\vartheta+r)^2}{4r}$, \cr
L_p \cdot s_p \cdot p^{-\min\{(\sqrt{r}-\sqrt{q})^2, \vartheta
+[(\sqrt{r(1-a^2)}-\sqrt{q})^+]^2 \}},
&\quad if $\dfrac{(\vartheta+r)^2}{4r}< q<r$, \cr
s_p \cdot\bigl(1+o(1)\bigr), &\quad if $q > r$.}
\]
\end{lemma}

The exponents on the right-hand side are minimized at $q =(\vartheta+
r)^2/(4r)$ if $r/\vartheta< [2 - \sqrt{1 - a^2}]/[\sqrt{1 - a^2} ( 1
- \sqrt{1 - a^2}) ]$, and at\vspace*{1pt}
$q=[2\vartheta+r(1-a^2)]^2/\allowbreak[4r(1-a^2)]$ if $r/\vartheta> [2 - \sqrt
{1 - a^2}]/[\sqrt{1 - a^2} ( 1 - \sqrt{1 - a^2}) ]$.
%
As a result, we have the following theorem, the proof of which is omitted.
%
%
\begin{theorem} \label{thmsubset}
Set the tuning parameter $\lambda^{\mathrm{ss}}_p = \sqrt{2 q \log p}$. Then
for all $q > 0$, the Hamming error of the subset selection satisfies
\begin{eqnarray*}
&&\frac{\hamm_p(\hb^{\mathrm{ss}}(\lambda_p^{\mathrm{ss}}))}{s_p}
\\
&&\qquad\geq
\cases{
L_p p^{-(\vartheta-r)^2/(4r)}, &\quad if $\displaystyle
\frac{r}{\vartheta} < \frac{2
- \sqrt{1 - a^2}}{\sqrt{1 - a^2} ( 1 - \sqrt{1 - a^2}) }$,
\vspace*{2pt}\cr
L_p p^{-{[2\vartheta+r(1-a^2)]^2}/({4r(1-a^2)}) +\vartheta}, &\quad if
$\displaystyle \frac{r}{\vartheta} > \frac{2 - \sqrt{1 - a^2}}{\sqrt{1 - a^2} (
1 - \sqrt{1 - a^2})
}$.}
\end{eqnarray*}
\end{theorem}

This gives\vspace*{1pt} the phase diagram in Figure \ref{figSSPhase}, where
$(\vartheta, r)$ satisfying
$r/\vartheta< [2 - \sqrt{1 - a^2}]/[\sqrt{1 - a^2} ( 1 - \sqrt{1 -
a^2}) ]$ defines the\vspace*{1pt} optimal region, and $(\vartheta, r)$ with
$r/\vartheta> [2 - \sqrt{1 - a^2}]/[\sqrt{1 - a^2} ( 1 - \sqrt{1 -
a^2}) ]$ defines the nonoptimal region.
Similar to the lasso, the subset selection is able to separate signal
singletons from the pure noise provided that the tuning parameter is
properly set.
But the subset selection is too harsh on signal pairs, triplets, etc.,
which costs its rate optimality.
In \cite{Thesis}, we further show that in certain regions of the phase
space, the lower bound in Theorem \ref{thmlasso} is tight.

\section{Simulations} \label{secSimul}

We have conducted a small-scale empirical study of the performance of
the UPS. The idea is to select a few interesting combinations of
$(\vartheta, \theta, \pi_p, \Omega)$ and study the behavior of the
UPS for finite $p$.
Fixing $(p, \pi_p, \Omega, \vartheta, \theta)$, let $n_p =
p^{\theta}$ and $\eps_p = p^{- \vartheta}$.
We investigate both the random design model and Stein's normal means model.

In the former, the experiment contains the following steps: (1)
Generate a \mbox{$p\times1$} vector $\beta$ by $\beta_j
\stackrel{\mathrm{i.i.d.}}{\sim} (1 - \eps_p) \nu_0 + \eps_p \pi_p$,
and an
$n_p\times1$ vector $z \sim N(0, I_{n_p})$. (2)~Generate an
$n_p\times p$ matrix $X$ the rows of which are samples from $N(0,
\frac{1}{n_p} \Omega)$; let $Y = X \beta+ z$. (3) Apply the UPS and
the lasso. For the lasso, we use the \textit{glmnet} package by Friedman
et al. \cite{Friedman2010} ($\Omega$ is assumed unknown in both
procedures).
(4) Repeat 1--3 for $100$ independent cycles, and calculate the average
Hamming distances.

In the latter, the settings are similar, except for (i) $n_p = p$, (ii)
$Y \sim N(\Omega^{1/2} \beta, I_p)$ in step 2 and (iii)
$\Omega$ is assumed as known in step 3 (otherwise valid inference is
impossible).
We include Stein's normal means model in the study for it is the
idealized version of the random design model.
\begin{Experiment}\label{exper1}
In this experiment, we use Stein's normal means
model to investigate the boundaries of the region of exact recovery by
the UPS and that by the lasso. Fixing $p = 10^4$ and $\Omega$ as the
tridiagonal matrix in~(\ref{tridiagonal2}) with $a = 0.45$, we let
$\vartheta$ range in $\{ 0.25, 0.5, 0.65\}$, and let \mbox{$\pi_p = \nu
_{\tau_p}$} with $\tau_p = \sqrt{2 r \log p}$, where $r$ is chosen
such that $\tau_p \in\{5, 6, \ldots, 12\}$.
For both procedures, we use the ideal threshold introduced in Sections
\ref{secUPS} and \ref{seclasso}, respectively.
That is, the tuning parameters of the UPS are set as $(t_p^*, \lambda
_p^{\mathrm{ups}}, u_p^{\mathrm{ups}}) = (\frac{\vartheta+ r}{2r}\tau_p, \sqrt{2
\vartheta\log(p)}, \tau_p)$, and the
tuning parameter of the lasso is set as $\lambda^{\mathrm{lasso}}_p = \max\{
\frac{\vartheta+ r}{2r}, (1 + \sqrt{(1 - a)/(1 + a)})^{-1} \}
\tau_p$.\vspace*{1pt}

%
\begin{table}
\caption{Hamming errors (Experiment \protect\ref{exper1}). UPS needs weaker signals for
exact recovery}
\label{tabexactrecbound}
\begin{tabular*}{\tablewidth}{@{\extracolsep{\fill}}lcd{3.2}d{2.2}
d{2.2}d{2.2}d{2.2}d{2.2}d{2.2}d{2.2}@{}}
\hline
& $\bolds{\tau_p}$ & \multicolumn{1}{c}{\textbf{5}}
& \multicolumn{1}{c}{\textbf{6}} & \multicolumn{1}{c}{\textbf{7}}
& \multicolumn{1}{c}{\textbf{8}} & \multicolumn{1}{c}{\textbf{9}}
& \multicolumn{1}{c}{\textbf{10}} & \multicolumn{1}{c}{\textbf{11}}
& \multicolumn{1}{c@{}}{\textbf{12}} \\
\hline
$\vartheta= 0.25$ & UPS & 49 & 11.1 & 1.79 & 0.26 &
0.02 & 0 & 0 & 0 \\
& lasso & 186.7 & 99.35 & 58.26 & 38.53 & 25.97 & 18.18 & 12.94 & 10.57
\\
[4pt]
$\vartheta=0.50$ & UPS &10.06 & 2.11 & 0.37 & 0.09 &
0 & 0 & 0 & 0 \\
& lasso & 16.36 & 5.11 & 1.47 & 0.51 & 0.28 & 0.33 & 0.26 & 0.09\\
[4pt]
$\vartheta= 0.65$ & UPS & 5.49 & 1.29 & 0.33 & 0.06
& 0 & 0 & 0 & 0\\
& lasso & 7.97 & 2.43 & 0.69 & 0.18 & 0.07 & 0.03 & 0.02 & 0.01\\
\hline
\end{tabular*}
\end{table}

The results are reported in Table \ref{tabexactrecbound}, where
the UPS outperforms consistently over the lasso, most prominently in
the case of $\vartheta= 0.25$. Also, for $\vartheta= 0.25, 0.5, $ or
$0.65$, the Hamming errors of the UPS start to fall below~$1$ when
$\tau_p$ exceeds $8, 7$ or $7$, respectively, but that of the lasso
won't fall below~$1$ until $\tau_p$ exceeds $12, 8$ or $ 7$,
respectively. In Section \ref{secIntro}, we show that the UPS yields
exact recovery when $\tau_p > (1 + \sqrt{1 - \vartheta}) \sqrt{2
\log p}$, where the right-hand side equals $(8.01, 7.32, 7.01)$ with
the current choices of $(p, \vartheta)$. The numerical results fit well
with the theoretic results.\vadjust{\goodbreak}
\end{Experiment}
\begin{Experiment}\label{exper2}
We use a random design model where $(p, \vartheta,
\theta) = (10^4,\allowbreak 0.65$, $0.91)$, and $\tau_p \in\{1, 2, \ldots, 7\}
$. The experiment contains three parts, 2a--2c. In 2a, we take $\Omega
$ to be the penta-diagonal matrix $\Omega(i,j) = 1\{ i = j \} + 0.4
\cdot1\{ |i - j| = 1\} + 0.1 \cdot1\{ |i - j| = 2\}$. Also,
for each $\tau_p$, we set $\pi_p$ as $\operatorname{Uniform}(\tau_p -
0.5, \tau_p + 0.5)$.
In 2b, we generate $\Omega$ in a way such that it has $4$ nonzero
off-diagonal elements on average in each row and each column, at
locations randomly chosen. Also, for each $\tau_p$, we take $\pi_p$
to be $\operatorname{Uniform}(\tau_p - 1, \tau_p + 1)$.
In 2c, we use a non-Gaussian design for $X$. In detail, first, we
generate an $n \times p$ matrix
$M$ the coordinates of which are i.i.d. samples from
$\operatorname{Uniform}(-\sqrt{3}, \sqrt{3})$. Second, we generate
$\Omega$ as in
2b. Last, we let $X = (1/\sqrt{n}) M \Omega^{1/2}$.
%
%
\begin{table}[b]
\caption{Ratios between Hamming errors and $p \eps_p$ (Experiment
\protect\ref{exper2}\textup{a}--\protect\ref{exper2}\textup{c}). Bold: UPS. Plain: lasso}
\label{tabex2}
\begin{tabular*}{\tablewidth}{@{\extracolsep{\fill}}lccccccc@{}}
\hline
$\bolds{\tau_p}$ & \multicolumn{1}{c}{\textbf{1}} & \multicolumn{1}{c}{\textbf{2}}
& \multicolumn{1}{c}{\textbf{3}} & \multicolumn{1}{c}{\textbf{4}}
& \multicolumn{1}{c}{\textbf{5}} & \multicolumn{1}{c}{\textbf{6}}
& \multicolumn{1}{c@{}}{\textbf{7}} \\
\hline
2a & \textbf{1.01} 1.02 & \textbf{0.96} 1.04 & \textbf{0.82} 0.97
& \textbf{0.51} 0.64 & \textbf{0.24} 0.28 & \textbf{0.09} 0.10
& \textbf{0.04} 0.04 \\
2b & \textbf{1.00} 1.00 & \textbf{0.98} 1.04 & \textbf{0.84} 0.96
& \textbf{0.55} 0.67 & \textbf{0.26} 0.32
& \textbf{0.10} 0.12& \textbf{0.05} 0.05\\
2c & \textbf{0.94} 0.95 & \textbf{0.90} 0.91 & \textbf{0.89} 0.95
& \textbf{0.48}
0.60 & \textbf{0.18} 0.27 & \textbf{0.05} 0.11 & \textbf{0.01} 0.03\\
\hline
\end{tabular*}
\end{table}
Also, for each $\tau_p$, we take $\pi_p$ to be the mixture of two
uniform distributions $\frac{1}{2} \operatorname{Uniform}(\tau_p - 0.5,
\tau_p+0.5) + \frac{1}{2} \operatorname{Uniform}(- \tau_p - 0.5, - \tau
_p + 0.5)$.
In all these experiments, the tuning parameters are set the same way as
in Experiment \ref{exper1}. The results are reported in Table \ref{tabex2},
suggesting that the UPS outperforms the lasso almost over the whole
range of $\tau_p$.
\end{Experiment}
%
%
\begin{Experiment}\label{exper3}
The goal of this experiment is twofold. First, we
investigate the sensitivity of the UPS and the lasso with respect to
their tuning parameters. Second, we investigate the refined UPS
introduced in Section~\ref{secrefinement}.
Fix $q > 0$. For the lasso, we take $\lambda^{\mathrm{lasso}}_p = \sqrt{2 q
\log(p)}$. For the UPS, set the $U$-step tuning parameter as $t_p^* =
\sqrt{2 q \log(p)}$ and let the $P$-step tuning parameters be
estimated as in
(\ref{UPStuningAdd}). Theorem \ref{thmadaptive} predicts that the UPS
performs well provided that $q \in(\max\{ \vartheta, \delta_0^2
(1+\eta)^2 r \}, (\vartheta+ r)^2/(4r))$, so both the lasso and
the UPS are driven by one tuning parameter $q$.
We now investigate how the choice of $q$ affects the performances of
the UPS and the lasso.
The experiment contains three sub-experiments 3a--3c.

%
\begin{figure}

\includegraphics{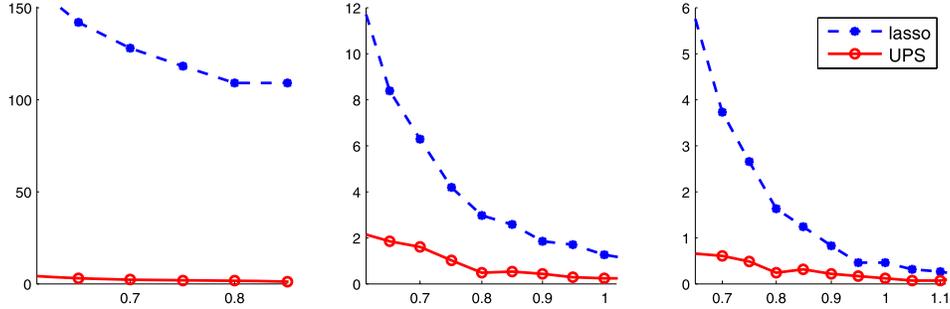}

\caption{Experiment \protect\ref{exper3}\textup{a}. $x$-axis: $q$. $y$-axis: Hamming error.
Left to right: $\vartheta= 0.2, 0.5, 0.65$.} \label{figEx3a2}
\end{figure}

In 3a, we use Stein's normal means model where $(p, r) = (10^4, 3)$,
\mbox{$\pi_p = \nu_{\tau_p}$} with $\tau_p = \sqrt{2 r \log p}$, $\Omega
$ is the penta-diagonal matrix satisfying $\Omega(i,j) = 1_{\{ i = j\}
} + 0.45 \cdot1_{\{ |i - j| = 1\}} + 0.05 \cdot1_{\{ |i - j| = 2\}}$,
and $\vartheta\in\{0.2, 0.5, 0.65\}$.
Note that when $\vartheta= 0.65$, $(\max\{ \vartheta, \delta_0^2
(1+\eta)^2 r \}, (\vartheta+ r)^2/(4r)) = (0.65, 1)$ (similarly
for other~$\vartheta$), so we let $q \in\{ 0.7, 0.8, \ldots, 1.1\}$.

In 3b, we use a random design model where $(p, r, \pi_p, \Omega, q)$
and the tuning parameters are the same as in 3a, but $\theta= 0.8$ and
$\vartheta\in\{0.5, 0.65\}$ (the case $\vartheta= 0.2$ is relatively
challenging in computation so is omitted).
We compare the lasso with the refined UPS where in each iteration, we
use the same tuning parameters as in~3a.

%
\begin{figure}[b]

\includegraphics{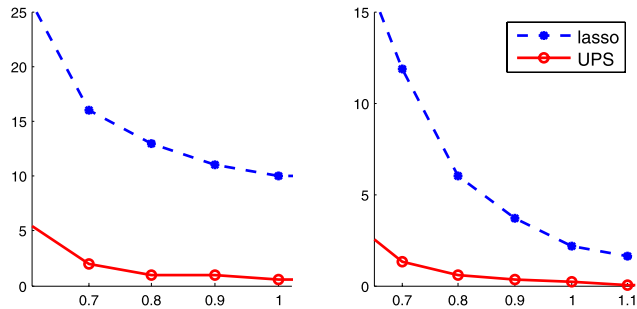}

\caption{Experiment \protect\ref{exper3}\textup{b}. $x$-axis: $q$. $y$-axis: Hamming error. Left:
$\vartheta= 0.5$. Right: $\vartheta= 0.65$. }
\label{figEx3a}
\end{figure}

In 3c, we use the same setup as in 3b, except that we fix $q = 1$ and
let $\tau_p$ range in $\{6, 6.5, \ldots, 9\}$.

The results of 3a--3c are reported in Figures \ref{figEx3a2}--\ref
{figEx3b}, correspondingly. These results suggest that, first, the UPS
consistently outperforms the lasso, and, second, the
UPS is relatively less sensitive to different choices of $q$.
\end{Experiment}
\begin{Experiment}\label{exper4}
In this experiment, we investigate the effect of
larger~$p$ and~$n$, respectively.
The experiment includes two sub-experiments, 4a and~4b.

In 4a, we use Stein's normal means model where $(\vartheta, r)=(0.5,
3)$, $\Omega$~as in Experiment \ref{exper2}c, $\pi_p=\nu_{\tau_p}$ with $\tau
_p = \sqrt{2 r \log p}$, and we let $p = 100 \times\{1, 10$, $
10^2,10^3, 10^4\}$. The lasso and the UPS are implemented as in
Experiment~\ref{exper3}a, where $q = 1$. The results are reported in the left part
of Table \ref{tabnpa}, where the second line displays the ratios
between the Hamming errors by the lasso and that by the UPS. Theoretic
results (Sections \ref{subseclasso} and \ref{seclasso}) predict that
for $(\vartheta, r)$ in the nonoptimal region of the lasso, such
ratios diverge as $p$ tends to $\infty$. The numerical results fit
well with the theory.

%
\begin{figure}

\includegraphics{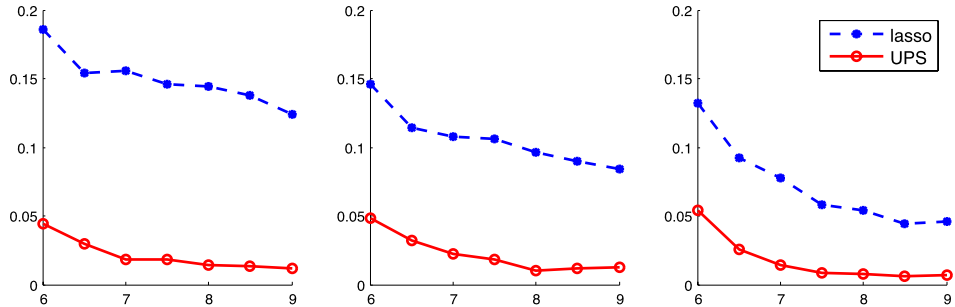}

\caption{Experiment \protect\ref{exper3}\textup{c}. The $x$-axis is $\tau_p$, and the $y$-axis is
the ratio between the Hamming error and $p\varepsilon_p$. Left to right:
$\vartheta=0.65, 0.5, 0.2$.}
\label{figEx3b}
\end{figure}

%
\begin{table}[b]
\caption{Left: ratios between the Hamming errors by the UPS and that by
the lasso (Experiment~\protect\ref{exper4}\textup{a}). Right: ratios between the
Hamming errors by the UPS for the random design model and that for
Stein's normal means model (Experiment~\protect\ref{exper4}\textup{b})}
\label{tabnpa}
\begin{tabular*}{\tablewidth}{@{\extracolsep{\fill}}lccccccccc@{}}
\hline
\multicolumn{5}{@{}c}{$\bolds{p}$} & \multicolumn{5}{c@{}}{$\bolds{n}$}\\[-4pt]
\multicolumn{5}{@{}c}{\hrulefill} & \multicolumn{5}{c@{}}{\hrulefill}\\
$\bolds{10^2}$ & $\bolds{10^3}$ & $\bolds{10^4}$ &
$\bolds{10^5}$ & $\bolds{10^6}$ & \textbf{300} & \textbf{900}
& \textbf{2,700} & \textbf{8,100} & \textbf{24,000}\\
\hline
2.43 & 5.81& 6.25& 8.80& 10.37& 479.25 & 54.04 & 12.66 & 1.08 &
1.01\\
\hline
\end{tabular*}
\end{table}

In 4b, we illustrate that in a random design model, if we fix $p$ and
let $n$ increase, then the
random design models get increasingly close to Stein's normal means model.
In detail,
we take a random design model where $(p, \vartheta,r)=(10^4, 0.5, 3)$,
$\Omega$ and $\pi_p$ as in Experiment \ref{exper2}c and $n_p = 300 \times\{1,
3, 3^2, 3^3, 3^4\}$.
We also take Stein's normal means model with the same $(p, \vartheta,
r, \Omega, \pi_p)$.
The performance of the UPS in both models is reported in the right part
of Table~\ref{tabnpa}, where the last line is the ratio between the
Hamming errors by the UPS for the random design model and that for the
Stein's normal means model.
The ratios effectively converge to $1$ as $n$ increases.
\end{Experiment}

%

\section*{Acknowledgments}

Jiashun Jin thanks Tony Cai, Emmanuel Candes, David Donoho, Stephen Fienberg,
Alan Friez, Robert Nowak, Runze Li, Larry Wasserman and Cun-Hui Zhang
for valuable pointers and discussion.

\begin{supplement}[id=suppA]
\stitle{Supplementary material for ``UPS delivers optimal phase diagram
in high-dimensional variable selection''}
\slink[doi]{10.1214/11-AOS947SUPP} 
\sdatatype{.pdf}
\sfilename{aos947\_supp.pdf}
\sdescription{Owing to space constraints, the technical proofs are
moved to a supplementary document~\cite{UPSSupp}.}
\end{supplement}

%

\printaddresses

\end{document}